\documentclass[10pt]{article}
\pdfoutput=1
\usepackage{geometry, amsmath, amsthm, latexsym, amssymb, graphicx, mathtools, tikz, enumerate, dsfont}
\usepackage{tikz-cd}
\usepackage[utf8]{inputenc}
\usepackage[T1]{fontenc}
\usepackage[shortlabels]{enumitem}
\setenumerate{topsep=0pt}

\usetikzlibrary{arrows.meta,positioning,calc}
\usepackage{tkz-graph}
\usepackage{graphicx}
\usepackage{csquotes}
\newcommand{\git}{\mathbin{
  \mathchoice{/\mkern-6mu/}
    {/\mkern-6mu/}
    {/\mkern-5mu/}
    {/\mkern-5mu/}}}
\usepackage{lmodern}
\usepackage{arydshln}
\usepackage{dynkin-diagrams}
\usepackage[shortlabels]{enumitem}
\usepackage{multirow}
\usepackage{mathdots}
\usepackage{hyperref}

\theoremstyle{definition}
\newtheorem{defi}{Definition}[section]
\newtheorem{ex}[defi]{Example}
\newtheorem{obs}[defi]{Remark}

\theoremstyle{plain}
\newtheorem{theorem}[defi]{Theorem}
\newtheorem{prop}[defi]{Proposition}
\newtheorem{lemma}[defi]{Lemma}
\newtheorem{cor}[defi]{Corollary}
\newtheorem*{idea*}{Idea}

\DeclareMathOperator{\SL}{SL}

\DeclareMathOperator{\PGL}{PGL}

\DeclareMathOperator{\So}{SO}
\DeclareMathOperator{\PSo}{PSO}

\DeclareMathOperator{\id}{id}

\DeclareMathOperator{\ord}{ord}

\DeclareMathOperator{\Aut}{Aut}
\DeclareMathOperator{\Isom}{Isom}

\DeclareMathOperator{\Hom}{Hom}

\DeclareMathOperator{\Ad}{Ad}
\DeclareMathOperator{\Sym}{Sym}
\DeclareMathOperator{\Spec}{Spec}
\DeclareMathOperator{\Ht}{ht}
\DeclareMathOperator{\ad}{ad}
\DeclareMathOperator{\diag}{diag}
\DeclareMathOperator{\Stab}{Stab}
\DeclareMathOperator{\Lie}{Lie}
\DeclareMathOperator{\End}{End}
\DeclareMathOperator{\vol}{vol}

\DeclareMathOperator{\cl}{cl}
\DeclareMathOperator{\Res}{Res}

\newcommand{\lp}{\left(}
\newcommand{\rp}{\right)}

\newcommand{\la}{\langle}
\newcommand{\ra}{\rangle}
\newcommand{\Z}{\mathbb{Z}}
\newcommand{\Q}{\mathbb{Q}}
\newcommand{\R}{\mathbb{R}}
\newcommand{\C}{\mathbb{C}}
\newcommand{\A}{\mathbb{A}}
\newcommand{\G}{\mathbb{G}}

\newcommand{\N}{\mathbb{N}}

\newcommand{\oO}{\mathcal{O}}

\newcommand{\pp}{\mathfrak{p}}
\newcommand{\Gg}{\mathfrak{g}}
\newcommand{\hh}{\mathfrak{h}}
\newcommand{\fz}{\mathfrak{z}}
\newcommand{\fc}{\mathfrak{c}}
\newcommand{\ft}{\mathfrak{t}}
\newcommand{\fsl}{\mathfrak{sl}}

\newcommand{\cC}{\mathcal{C}}

\newcommand{\cR}{\mathcal{R}}
\newcommand{\cL}{\mathcal{L}}

\newcommand{\cW}{\mathcal{W}}
\newcommand{\cS}{\mathcal{S}}
\newcommand{\cF}{\mathcal{F}}

\newcommand{\cU}{\mathcal{U}}
\newcommand{\cB}{\mathcal{B}}
\newcommand{\cE}{\mathcal{E}}
\newcommand{\eps}{\varepsilon}
\newcommand{\ul}{\underline}
\newcommand{\ol}{\overline}

\usepackage[
backend=biber,
style=alphabetic,
sorting=nyt
]{biblatex}
\title{Geometry-of-numbers over number fields and the density of 
ADE families of curves having squarefree discriminant}
\author{Martí Oller}

\begin{document}

\maketitle
\begin{abstract}
For families of curves arising from a Dynkin diagram of type ADE, we show
that the density of such curves having squarefree discriminant is equal
to the product of local densities. We do so using the framework of Thorne 
and Laga's PhD theses and geometry-of-numbers techniques developed
by Bhargava, here expanded over number fields.
\end{abstract}
\tableofcontents
\section{Introduction}

In this paper, we aim to determine the density of curves in certain families
that have squarefree discriminant.
Over $\Q$, this was done in \cite{Oller}, and in this article we will
generalise the methods in said paper to a general number field $F$. For
completeness and convenience, this article will be self-contained, in the
sense that it will develop the results in \cite{Oller} from scratch.

We will generalise the methods of the article \cite{BSWsquarefree} by 
Bhargava, Shankar and Wang, in which they compute the density of monic integral
polynomials of a given degree that have squarefree discriminant. 
In their situation, they relate polynomials with discriminant divisible
by the square of a large prime with orbits of the representation $(G,V)$,
$G = \So_n$ acts on the space $V$ of $n\times n$ symmetric matrices. 
The key observation that motivates our results is that the representation 
$(G,V)$ arises as a particular case of the more 
general families of representations studied in \cite{ThorneThesis}.
Given a simply laced Dynkin diagram, Thorne used Vinberg theory to associate to it a 
family of curves and a coregular representation $(G,V)$, where the rational
orbits of the representation are related to the arithmetic of the curves
in the family. This parametrisation has been used
to study the size of $2$-Selmer groups of the Jacobians of these curves, 
see \cite{BG,SW,ShankarD2n+1,ThorneE6,RTE78,LagaE6} for some particular 
cases. Later, Laga unified, reproved and extended all these results in \cite{LagaThesis} 
in a uniform way.

Our aim is to compute the density of curves having squarefree discriminant
in these families of $ADE$ curves. We will do so by reinterpreting the
methods in \cite{BSWsquarefree} in the language of \cite{ThorneThesis}
and \cite{LagaThesis}. Given that we prove our results over a general number field $F$, this presents
an additional challenge, in the sense that most of the literature on the
required geometry-of-numbers works over $\Q$, and the translation to the
number field case is not necessarily immediate. Taking \cite{BSWprehom}
and \cite{BSWcoregular} as a point of reference, we develop the techniques
that we need in geometry-of-numbers over a number field.

Let $\mathcal{D}$ be a Dynkin diagram of type $A,D,E$, and let $F$ be a number field. 
In Section \ref{subs:Vinberg}, we will construct a representation $(G,V)$ associated to $\mathcal{D}$.
Let $B$ denote the Geometric Invariant Theory (GIT)
quotient $V \git G := \Spec F[V]^G$, which is isomorphic to
an affine space: $B = \Spec F[p_{d_1},\dots,p_{d_r}]$.
The group $\G_m$ acts on $B$ by $\lambda \cdot p_{d_i} = \lambda^{d_i}p_{d_i}$.
We want to define a sensible notion of \emph{height} for elements of $b \in B(F)$.
In Section \ref{subs:heights}, for $b \in B(F)$ we will define
\[
\text{ht}(b) := (NI_b) \prod_{v \in M_{\infty}} \sup\left(|p_{d_1}(b)|_{v}^{1/d_1},\dots,|p_{d_r}(b)|_{v}^{1/d_r}\right),
\]
where $M_{\infty}$ denotes the set of archimedean places of $F$ and $I_b$
is the ideal $I_b = \{a \in F \mid a^{d_i}p_{d_i}(b) \in \oO_F, \, \forall i\}$.
This height is $\G_m(F)$-invariant by the product formula.
Further, we can see that the number of elements of $\G_m(F)\backslash B(F)$
having bounded height is finite; see \cite[Theorem A]{Deng} for a more precise count.
Throughout the paper, we will fix a subset $\Sigma \subset B(\oO_F)$ which is a fundamental
domain for the action of $\G_m(F)$ on $B(F)$: we will construct it in
Section \ref{subs:heights}.

In Section \ref{subs:Transverse}, we will construct a family of curves $C \to B$
associated to the given Dynkin diagram.
Denote by $C_b$ the preimage of a given $b\in B$ under the map $C \to B$;
it will be a curve of the form given by Table \ref{table:curves}. The main
result of this paper concerns the density of squarefree values of the
discriminant $\Delta(C_b)$ of the curve: a definition for the
discriminant of a plane curve can be found in \cite[\S 2]{Sutherland2019}, for instance.
We remark that in our definition of discriminant, we assume that it is
an polynomial in multiple variables and coefficients in $\oO_F$, normalised so that
the coefficients have common divisors (for instance, the usual
discriminant for elliptic curves contains a factor of $16$: we omit it
in our case). Alternatively, by Theorem \ref{theo:basics} our main result
can be interpreted as a statement about the invariant polynomial $\Z[\hh]^H$
for any simple adjoint algebraic group $H$ of type $ADE$.

For a prime ideal $\pp$ of $\oO_F$, we will denote the completion of $\oO_{F}$
with respect to $\pp$ by $\oO_{\pp}$. We will further set $F_{\pp}$ to
be the field of fractions of $\oO_{\pp}$ and $k_{\pp}$ to be the corresponding
residue field. Now, consider the set
\[
B'(\oO_{\pp}) = \{b \in B(\oO_{\pp}) \mid v_{\pp}(p_{d_i}(b)) < d_i \text{ for some }i\}.
\]
Every element of $B(F_{\pp})$ is $\G_m(F_{\pp})$-conjugate to an element
of $B'(\oO_{\pp})$, and this element is unique up to the $\G_m(\oO_{\pp})$-action.
For a given $b \in B(F_{\pp})$, we say that $\Delta(b)$ is \emph{squarefree at $\pp$} if, for
any element $b' \in B'(\oO_{\pp})$ which is $\G_m(F_{\pp})$-conjugate to $b$,
the discriminant $\Delta(b')$ is squarefree as an element of $\oO_{\pp}$.
For $b \in B(F)$, we say that $\Delta(b)$ is \emph{squarefree} if $\Delta(b)$ is squarefree at $\pp$
for all finite primes $\pp$. Note that the property of ``being squarefree''
does not change with the action of $\G_m(F)$.

Our result is related to the $\pp$-adic density of these squarefree values: 
we will denote by $\rho(\mathcal{D}_{\pp})$ the local density at $\pp$ of curves in 
the family $C \to B$ having discriminant indivisible by $\pp^2$ in $F_{\pp}$;
this is obtained by taking all the (finitely many) elements in $b \in B(F_{\pp}/\pp^2F_{\pp})$
and counting the proportion of them that have non-zero discriminant in
$F_{\pp}/\pp^2F_{\pp}$. We note that under our assumptions on the discriminant none 
of the local densities vanish; this can be checked with a case-by-case computation.
\begin{theorem}
\label{theo:main}
We have
\[
\lim_{X \to \infty} \frac{\#\{b \in \G_m(F)\backslash B(F) \mid \Delta(b) \text{ is squarefree, } \Ht(b) < X\}}{\#\{b \in \G_m(F)\backslash B(F) \mid \Ht(b) < X\}} = \prod_{\pp} \rho(\mathcal{D}_{\pp}).
\]
\end{theorem}
To prove Theorem \ref{theo:main}, we need to obtain a tail estimate to show that
``not too many'' $b \in B(\oO_F)$ have discriminant divisible by $I^2$ for squarefree ideals $I$ of large norm.
Let $\pp$ be a prime ideal. A key observation in \cite{BSWsquarefree} is to separate those $b$
with $\pp^2$ dividing $\Delta(b)$ in two separate cases:
\begin{enumerate}
\item If $\pp^2 | \Delta(b+pc)$ for all $p \in \pp$ and $c \in B(\oO_F)$, we say $\pp^2$ \emph{strongly divides}
$\Delta(b)$ (in other words, $\pp^2$ divides $\Delta(b)$ for ``mod $\pp$ reasons'').
\item If there exists $p \in \pp$ and $c \in B(\oO_F)$ such that $\pp^2 \nmid \Delta(b+pc)$,
we say $\pp^2$ \emph{weakly divides} $\Delta(b)$ (in other words, $\pp^2$ 
divides $\Delta(b)$ for ``mod $\pp^2$ reasons'').
\end{enumerate}
Analogously, for a squarefree ideal $I \subset \oO_F$, we will say that
$I^2$ strongly (resp. weakly) divides $\Delta(b)$ if every prime ideal
$\pp \mid I$ strongly (resp. weakly) divides $\Delta(b)$.
We will let $\cW_{I}^{(1)},\cW_{I}^{(2)}$ denote the set 
of $b \in B(\oO_{F})$ whose discriminant is strongly (resp. weakly) divisible 
by $I^2$. We want to prove tail estimates for $\cW_{I}^{(1)},\cW_{I}^{(2)}$ separately.
Our argument for the weakly divisible case will require us to avoid finitely
many primes: more precisely, in Section \ref{section:construct} we will define an element
$N_{bad} \in \oO_F$ which will be divisible by all these ``bad primes''.
\begin{theorem}
\label{theorem:tail}
There exists a constant $\delta > 0$ such that for any positive real number $M$ we have:
\begin{align*}
\sum_{\substack{I \text{ squarefree} \\ NI > M}} \# \{b \in \G_m(F)\backslash \cW_{I}^{(1)} \mid \Ht(b) < X\} = O_{\eps}\left(\frac{X^{\dim V + \eps}}{M}\right) + O_{\eps}\left(X^{\dim V - 1+\eps}\right),\\
\sum_{\substack{I \text{ squarefree} \\ NI > M \\ (I,N_{bad}) = 1}} \#\{b \in \G_m(F)\backslash \cW_{I}^{(2)} \mid \Ht(b) < X\} = O_{\eps}\left(\frac{X^{\dim V+\eps}}{M}\right) + O\left(X^{\dim V - \delta}\right).
\end{align*}
The implied constants are independent of $X$ and $M$.
\end{theorem}
The strongly divisible case will
follow from the use of the Ekedahl sieve: see Section \ref{subs:tail} for
a discussion. Hence, we will spend most of the paper dealing with the
weakly divisible case.

We start in Section \ref{section:Prelim}, where we develop the necessary
background and introducing our objects of interest, most importantly the representation
$(G,V)$ coming from Vinberg theory and the associated family of curves
$C \to B$. The main step in the proof of Theorem \ref{theorem:tail} is
done in Section \ref{section:construct}, where given an element $b \in \cW_I^{(2)}$,
we construct a special integral orbit in $V$, whose elements have invariant
$b$. We additionally consider a distinguished subspace $W_0 \subset V$, 
and we define a $Q$-invariant for the elements of $W_0$. Then, we will
see that the elements in the constructed orbit have large $Q$-invariant
when they intersect $W_0$ (which happens always except for a negligible
amount of times by cutting-off-the-cusp arguments). This construction
is the analogue of \cite[\S 2.2, \S 3.2]{BSWsquarefree}; we give a more detailed
comparison at the end of Section \ref{section:construct}.

In Section \ref{section:Reduction}, we set up the main tools that we will
require in reduction theory. This includes an extended discussion about
heights, as well as a construction of a suitable ``box-shaped'' fundamental
domain for the action of an arithmetic subgroup $\Gamma$ of $G(F)$ over
$G(F_{\infty})$, where $F_{\infty} = \prod_{v \in M_{\infty}} F_v$. Then, in Section \ref{section:Counting}, we compute a
precise asymptotic for the number of reducible orbits of bounded height.
This is a necessary step in the argument, since the proof of the main
results requires a power-saving asymptotic on the elements with big
stabiliser (Proposition \ref{prop:bs}), which will require the results
in Section \ref{section:Counting}. We remark that part of the argument
relies on extensive case-by-case computations: some of them are carried
out in Section \ref{subs:cases}, while some others take place implicitly
in the proof of Proposition \ref{prop:cusp}.
Finally, in Section \ref{section:final} we conclude the proof of the main
results. In Section \ref{subs:tail} we prove Theorem \ref{theorem:tail},
and in Section \ref{subs:sieve} we deduce Theorem \ref{theo:main} using
a squarefree sieve. 

\paragraph*{Acknowledgements.}
This paper was written while the author was a PhD student under the supervision
of Jack Thorne. I would like to thank him for providing many useful
suggestions, guidance and encouragement during the process, and for revising
an early version of this manuscript. I also wish to thank Jef Laga for
his helpful comments. I would also like to thank Manjul Bhargava, Arul
Shankar and Xiaoheng Wang for sharing the manuscript \cite{BSWcoregular}
with me.

The project 
that gave rise to these results received the support of a fellowship 
from ``la Caixa'' Foundation (ID 100010434). The fellowship code is LCF/BQ/EU21/11890111.
The author wishes to thank them, as well as the Cambridge Trust and the DPMMS,
for their support.

\subsection{Notation}
We recap the most important bits of notation in this section. Most of it
has already been introduced or will be introduced in the future, but is
included here for the convenience of the reader.

Throughout, we will work with a fixed number field $F$. We will denote
its ring of integers by $\oO_F$. For a finite prime $\pp$ of $F$, we will
denote by $\oO_{\pp}$ the completion of $\oO_F$ with respect to $\pp$,
$F_{\pp}$ its field of fractions and $k_{\pp}$ its residue field. We will
denote the set of infinite places of $F$ by $M_{\infty}$, and for any
$v \in M_{\infty}$, we will denote by $F_v$ the completion of $F$ with
respect to $v$. For an ideal $I \subset \oO_F$, we will denote by the norm
of the ideal by $NI = \#(\oO_F/I)$.

We will also denote $F_{\infty} = \prod_{v \in M_{\infty}} F_v$, and for
$x = (x_v)_v \in F_{\infty}$ we will denote
\[
|x| := \prod_{v}|x_v|_v,
\]
where $|x_v|_v$ denotes the norm in $F_v$ given by $|x_v|_v = N_{\C/F_v}(x_v)$.
For $x \in F$, we will denote by $|x|$ the norm of $F$ as an element of
$F_{\infty}$.

Given a split semisimple group $H$, we will consider a natural representation
$(G,V)$, where $G$ is a suitable subgroup of $H$. Inside $G$, we will fix
a split torus $T$, and a Borel subgroup $P$ containing $T$, corresponding
to the negative roots of $\Phi(G,T)$. We will also fix 
$B := V \git G = \Spec F[p_{d_1},\dots,p_{d_r}]$, where $p_{d_i}$ are
polynomials of degree $i$, on which $\G_m(F)$ acts upon by $\lambda \cdot p_{d_i} = \lambda^{d_i}p_{d_i}$.
We will also fix a Kostant section $\kappa \cdot B \to V$, which will
be a section of the invariant map $\pi \colon V \to B$.

We will also fix the unipotent radical $N$ of the Borel $P$, a maximal
compact subgroup $K$ of $G(F_{\infty})$ and the subgroup $A$ of $T(F_{\infty})$
consisting of $t = (t_1,\dots,t_k)$ such that $(t_i)_v \in \R_{> 0}$ for
every $i$ and every place $v \in M_{\infty}$.

\section{Preliminaries}
\label{section:Prelim}

In this section, we introduce our representation $(G,V)$ of interest,
together with some of its basic properties. We do so mostly following 
\cite[\S 2]{ThorneThesis} and \cite[\S 3]{LagaThesis}.

\subsection{Vinberg representations}
\label{subs:Vinberg}
Let $F$ be a number field, let $H_{\Q}$ be a split adjoint simple group of type $A,D,E$ over $\Q$,
and consider $H$ to be the base change of $H_{\Q}$ to $F$. We assume
$H$ is equipped with a pinning $(T_H,P_H,\{X_{\alpha}\})$, meaning:
\begin{itemize}
\item $T_H \subset H$ is a split maximal torus defined over $\Q$ (determining a root system $\Phi_{H}$).
\item $P_H \subset H$ is a Borel subgroup containing $T_H$ (determining a root basis $S_{H} \subset \Phi_{H}$).
\item $X_{\alpha}$ is a generator for $\hh_{\alpha}$ for each $\alpha \in S_{H}$.
\end{itemize}
Let $W = N_{H}(T_H)/T_H$ be the Weyl group of $\Phi_{H}$, and let $\mathcal{D}$ be
the Dynkin diagram of $H$. Then, we have the following exact sequences:
\begin{equation}
\label{eq:pin}
\begin{tikzcd}
0 \arrow{r} & H \arrow{r} & \Aut(H) \arrow{r} & \Aut(\mathcal{D}) \arrow{r} & 0
\end{tikzcd}
\end{equation}
\begin{equation}
\label{eq:ex2}
\begin{tikzcd}
0 \arrow{r} & W \arrow{r} & \Aut(\Phi_{H}) \arrow{r} & \Aut(\mathcal{D}) \arrow{r} & 0
\end{tikzcd}
\end{equation}
The subgroup $\Aut(H,T_H,P_H,\{X_{\alpha}\}) \subset \Aut(H)$ of automorphisms
of $H$ preserving the pinning determines a
splitting of \eqref{eq:pin}. Then, we can define $\vartheta \in \Aut(H)$
as the unique element in $(T_H,P_H,\{X_{\alpha}\})$ such that its image
in $\Aut(\mathcal{D})$ under \eqref{eq:pin} coincides with the image of $-1 \in \Aut(\Phi_{H})$
under \eqref{eq:ex2}. Writing $\check{\rho}$ for the sum of fundamental
coweights with respect to $S_{H}$, we define
\[
\theta := \vartheta \circ \Ad(\check{\rho}(-1)) = \Ad(\check{\rho}(-1)) \circ \vartheta.
\]
The map $\theta$ defines an involution of $H$, and so $d\theta$ defines an
involution of the Lie algebra $\hh$. By considering $\pm 1$ eigenspaces, we obtain a $\Z/2\Z$-grading
\[
\hh = \hh(0) \oplus \hh(1),
\]
where $[\hh(i),\hh(j)] \subset \hh(i+j)$. We define $G = (H^{\theta})^{\circ}$
and $V = \hh(1)$, which means that $V$ is a representation of $G$ by restriction
of the adjoint representation. Moreover, we have $\text{Lie}(G) = \hh(0)$.

We have the following basic result \cite[Theorem 1.1]{Panyushev} on the
GIT quotient $B := V \git G = \Spec F[V]^G$.

\begin{theorem}
\label{theo:basics}
Let $\mathfrak{c} \subset V$ be a Cartan subspace. Then, $\mathfrak{c}$
is a Cartan subalgebra of $\hh$, and the map $N_G(\mathfrak{c}) \to W_{\fc} := N_H(\fc)/Z_H(\fc)$
is surjective. Therefore, the canonical inclusions $\fc \subset V \subset \hh$
induce isomorphisms
\[
\fc \git W_{\fc} \cong V \git G \cong \hh \git H.
\]
In particular, all these quotients are isomorphic to a finite-dimensional affine space.
\end{theorem}

For any field $k$ of characteristic zero, we can define 
the \emph{discriminant polynomial} $\Delta \in k[\hh]^H$ as the image of 
$\prod_{\alpha \in \Phi_T} \alpha$ under the isomorphism $k[\mathfrak{t}]^{W} \xrightarrow{\sim} k[\hh]^H$.
The discriminant can also be regarded as a polynomial in $k[B]$ through the isomorphism
$k[\hh]^H \cong k[V]^G = k[B]$.
We can relate the discriminant to one-parameter subgroups, which we now introduce. If
 $\lambda \colon \G_m \to G_{k}$ is a homomorphism, there
exists a decomposition $V = \sum_{i \in \Z} V_i$, where $V_i := \{v \in V(k) \mid \lambda(t)v = t^iv, \; \forall t \in \G_m(k)\}$.
Every vector $v \in V(k)$ can be written as $v = \sum v_i$, where $v_i \in V_i$;
we call the integers $i$ with $v_i \neq 0$ the \emph{weights} of $v$.
Finally, we recall that an element $v \in \hh$ is \emph{regular} if
its centraliser has minimal dimension.
\begin{prop}
\label{prop:H-M}
Let $k/\Q$ be a field, and let $v \in V(k)$. The following are equivalent:
\begin{enumerate}
\item $v$ is regular semisimple.
\item $\Delta(v) \neq 0$.
\item For every non-trivial homomorphism $\lambda \colon \G_m \to G_{k^s}$,
$v$ has a positive weight with respect to $\lambda$.
\end{enumerate}
\end{prop}
\begin{proof}
The reasoning is the same as in \cite[Corollary 2.4]{RTE78}.
\end{proof}

We remark that the Vinberg representation $(G,V)$ can be identified explicitly. For the
reader's convenience, we reproduce the explicit description written in \cite[\S 3.2]{LagaThesis}
in Table \ref{table:explicit}. We refer the reader to \emph{loc. cit.} for the
precise meaning of some of these symbols.

\begin{table}[h]
\begin{center}
\begin{tabular}{l|c|c}
Type & $G$ & $V$ \\
\hline
$A_{2n}$ & $\So_{2n+1}$ & $\Sym^2(2n+1)_0$ \\
$A_{2n+1}$ & $\PSo_{2n+2}$ & $\Sym^2(2n+2)_0$ \\
$D_{2n}$ ($n \geq 2$) & $\So_{2n} \times \So_{2n} / \Delta(\mu_2)$ & $2n \boxtimes 2n$ \\
$D_{2n+1}$ ($n \geq 2$) &  $\So_{2n+1} \times \So_{2n+1}$ & $(2n+1)\boxtimes (2n+1)$ \\
$E_{6}$ & PSp$_8$ & $\wedge_{0}^4 8$ \\
$E_{7}$ & $\SL_{8}/\mu_4$ & $\wedge^4 8$ \\
$E_{8}$ & Spin$_{16}/\mu_2$& half spin\\
\end{tabular}
\caption{Explicit description of each representation}
\label{table:explicit}
\end{center}
\end{table}

\subsection{Restricted roots}
\label{section:ResRoots}

In the previous section we considered the root system $\Phi_H$ of $H$,
but we will also need to work with the root system of $G$. By \cite[Lemma 5.1]{Richardson},
we have that $T := T_H^{\theta}$ is a split maximal torus of $G$. We will
compare the root systems $\Phi_H$ and $\Phi_G = \Phi(G,T)$ in a similar
fashion to \cite[\S 2.3]{ThorneE6}.

Write $\Phi_H/\vartheta$ for the orbits of $\vartheta$ on $\Phi_H$, where
$\vartheta$ is the pinned automorphism defined in the previous section.
\begin{lemma}
\begin{enumerate}
\item The map $X^{\ast}(T_H) \to X^{\ast}(T)$ is surjective, and
the group $G$ is adjoint. In particular, $X^{\ast}(T)$ is spanned
by $\Phi_G$.
\item Let $\alpha, \beta \in \Phi_H$. Then, the image of $\alpha$ in $X^{\ast}(T)$
is non-zero, and $\alpha,\beta$ have the same image if and only if either
$\alpha = \beta$ or $\alpha = \vartheta(\beta)$.
\end{enumerate}
\end{lemma}
\begin{proof}
The fixed group $T$ is connected and contains regular elements
of $T_H$ by \cite[Lemma 3.1]{ReederTorsion}. The group $G$ has trivial centre
by \cite[\S 3.8]{ReederTorsion}. For the second part, see \cite[\S 3.3]{ReederTorsion}.
\end{proof}

Hence, we can identify $\Phi_H/\vartheta$ with its image in $X^{\ast}(T)$. We
note that $\vartheta = 1$ if and only if $-1$ is an element of the Weyl
group $W(H,T_H)$; in this case $\Phi_H/\vartheta$ coincides with $\Phi_H$.

We can write the following decomposition:
\[
\hh = \ft \oplus \bigoplus_{a \in \Phi_H/\vartheta} \hh_a,
\]
with $\ft = \ft^{\theta} \oplus V_0$ and $\hh_a = \Gg_a \oplus V_a$, according
to the $\theta$-grading. We have a decomposition
\[
V = V_0 \oplus \bigoplus_{a \in \Phi_H/\vartheta} V_a.
\]
For a given $a \in \Phi_H/\vartheta$ there are three cases to distinguish, 
according to the value of $s = (-1)^{\la \alpha,\check{\rho} \ra}$:
\begin{enumerate}
\item $a = \{\alpha\}$ and $s = 1$. Then, $V_a = 0$ and $\Gg_{\alpha}$ is spanned by $X_{\alpha}$.
\item $a = \{\alpha\}$ and $s = -1$. Then, $V_a$ is spanned by $X_{\alpha}$ and $\Gg_{\alpha} = 0$.
\item $a = \{\alpha,\vartheta(\alpha)\}$, with $\alpha \neq \vartheta(\alpha)$. 
Then, $V_a$ is spanned by $X_{\alpha}-sX_{\vartheta(\alpha)}$ and $\Gg_{\alpha}$ is spanned by $X_{\alpha}+sX_{\vartheta(\alpha)}$.
\end{enumerate}

We note that $\vartheta$ preserves the height of a root $\alpha$
with respect to the basis $S_H$. Therefore, it will make sense to define
the \emph{height} of a root $a \in \Phi/\vartheta$ as the height
of any element in $\vartheta^{-1}(a)$.

\subsection{Transverse slices over $V \git G$}
\label{subs:Transverse}

In this section, we present some remarkable properties of the map
$\pi \colon V \to B$, where we recall that $B := V \git G$ is the GIT quotient.

\begin{defi}
An \emph{$\fsl_2$-triple} of $\hh$ is a triple $(e,h,f)$ of non-zero elements
of $\hh$ satisfying
\[
[h,e] = 2e,\quad [h,f] = -2f, \quad [e,f] = h.
\]
Moreover, we say this $\fsl_2$-triple is \emph{normal} if $e,f \in \hh(1)$ and $h \in \hh(0)$.
\end{defi}
\begin{theorem}[Graded Jacobson-Morozov]
\label{theo:J-M}
Every non-zero nilpotent element $e \in \hh(1)$ is contained in a normal $\fsl_2$-triple.
If $e$ is also regular, then it is contained in a unique normal $\fsl_2$-triple.
\end{theorem}
\begin{proof}
The first part of the statement is \cite[Lemma 2.17]{ThorneThesis}, and
the second part follows from \cite[Lemma 2.14]{ThorneThesis}.
\end{proof}
\begin{defi}
Let $r$ be the rank of $\hh$.
We say an element $x \in \hh$ is \emph{subregular} if $\dim \fz_{\hh}(x) = r + 2$.
\end{defi}
Subregular elements in $V$ exist by \cite[Proposition 2.27]{ThorneThesis}.
Let $e \in V$ be such an element, and fix a normal $\fsl_2$-triple $(e,h,f)$
using Theorem \ref{theo:J-M}.
Let $C = e + \fz_{V}(f)$, and consider the natural morphism $\varphi \colon C \to B$.

\begin{theorem}
\label{theo:curves}
\begin{enumerate}
\item The geometric fibres of $\varphi$ are reduced connected curves.
For $b \in B(F)$, the corresponding curve $C_b$ is smooth if and only if
$\Delta(b) \neq 0$.
\item The central fibre $\varphi^{-1}(0)$ has a unique singular point which is a simple singularity
of type $A_n,D_n,E_n$, coinciding with the type of $H$. 
\item We can choose coordinates $p_{d_1},\dots,p_{d_r}$ in $B$, with $p_{d_i}$ being homogeneous
of degree $d_i$, and coordinates $(x,y,p_{d_1},\dots,p_{d_r})$ on $C$ such
that $C \to B$ is given by Table \ref{table:curves}.
\end{enumerate}
\end{theorem}
\begin{proof}
See \cite[Theorem 3.8]{ThorneThesis}.
\end{proof}

\begin{table}
\begin{center}
\begin{tabular}{l|l|c}
Type & Curve & \# Marked points\\
\hline
$A_{2n}$ & $y^2 = x^{2n+1} + p_2x^{2n-1} + \dots + p_{2n+1}$ & $1$\\
$A_{2n+1}$ & $y^2 = x^{2n+2} + p_2x^{2n} + \dots + p_{2n+2}$ & $2$ \\
$D_{2n}$ ($n \geq 2$) & $y(xy+p_{2n}) = x^{2n-1} + p_2x^{2n-2} + \dots + p_{4n-2}$ & $3$\\
$D_{2n+1}$ ($n \geq 2$) & $y(xy+p_{2n+1}) = x^{2n} + p_2x^{2n-1} + \dots + p_{4n}$ & $2$ \\
$E_{6}$ & $y^3 = x^4 + (p_2x^2 + p_5x + p_8)y + (p_6x^2 + p_9x+p_{12})$ & $1$\\
$E_{7}$ & $y^3 = x^3y + p_{10}x^2 + x(p_2y^2 + p_8y + p_{14}) + p_6y^2 + p_{12}y + p_{18}$ & $2$\\
$E_{8}$ & $y^3 = x^5 + (p_2x^3 + p_8 x^2 + p_{14}x + p_{20})y + (p_{12}x^3 + p_{18}x^2 + p_{24}x + p_{30})$ & $1$\\
\end{tabular}
\caption{Families of curves}
\label{table:curves}
\end{center}
\end{table}

Our choice of pinning in Section \ref{subs:Vinberg} determines a natural choice
of a regular nilpotent element, namely $e_{0} = \sum_{\alpha \in S_H} X_{\alpha} \in V(F)$.
Let $(e_0,h_0,f_0)$ be its associated normal $\fsl_2$-triple by Theorem \ref{theo:J-M}. We define the
affine linear subspace $\kappa_{e_0} := (e_{0} + \fz_{\hh}(f_{0})) \cap V$ as the
\emph{Kostant section} associated to $e_{0}$. Whenever $e_{0}$ is understood, we will just
denote the Kostant section by $\kappa$.

\begin{theorem}
The composition $\kappa \hookrightarrow V \to B$ is an isomorphism, and
every element of $\kappa$ is regular.
\end{theorem}
\begin{proof}
See \cite[Lemma 3.5]{ThorneThesis}.
\end{proof}

\begin{defi}
\label{defi:reducible}
Let $L/F$ be a field and let $v \in V(L)$. We say $v$ is \emph{$L$-reducible}
if $\Delta(v) = 0$ or if $v$ is $G(L)$-conjugate to some Kostant section,
and \emph{$L$-irreducible} otherwise.
\end{defi}
We will typically refer to $F$-(ir)reducible elements simply as (ir)reducible.
We note that if $L$ is algebraically closed, then all elements of $V$
are reducible, by \cite[Proposition 2.11]{LagaThesis}.

\subsection{Integral structures}
\label{subs:integral}
So far, we have considered our objects of interest over the fields $\Q$ and $F$, but for our
purposes it will be crucial to define integral structures for $G$ and $V$.

We start by considering structures over $\Z$. The structure of $G$ over $\Z$ comes from the general classification of
split reductive groups over any non-empty scheme $S$: namely, every root
datum is isomorphic to the root datum of a split reductive $S$-group
(see \cite[Theorem 6.1.16]{ConradRed}). By considering the root datum
$\Phi(G,T)$ studied in Section \ref{section:ResRoots} and the
scheme $S = \Spec \Z$, we get a split reductive group $\ul{G}_{\Z}$
defined over $\Z$, such that its base change to $\Q$ coincides with $G$.
By \cite[Lemma 5.1]{Richardson}, we know that $T$ is a maximal split torus
of $G$, and that $P := P_H^{\theta}$ is a Borel subgroup of $G$ containing
$T$.
We also get integral structures for $\ul{T}_{\Z}$ and $\ul{P}_{\Z}$
inside of $\ul{G}_{\Z}$. We get $\oO_F$-structures by base-changing to
$\oO_F$: set $\ul{G} := \ul{G}_{\Z} \otimes \oO_F$ and analogously
$\ul{T} := \ul{T}_{\Z} \otimes\oO_F$, $\ul{P} := \ul{P}_{\Z} \otimes \oO_F$.

For \emph{any} linear algebraic group $G$ defined over $F$, we recall that its class group
is
\[
\text{cl}(G) = (\prod_{\pp \notin M_{\infty}} G(\oO_{\pp}))\backslash G(\A_{F,f}) / G(F).
\]
\begin{prop}
\label{prop:clNum}
Every linear algebraic group $G$ has finite class group.
\end{prop}
\begin{proof}
See \cite[Theorem 5.1]{BorelCG}.
\end{proof}

To obtain a $\Z$-structure for $V_{\Q}$, we consider $\hh_{\Q}$ as a semisimple
$G_{\Q}$-module over $\Q$ via the restriction of the adjoint representation.
This $G_{\Q}$-module splits into a sum of simple $G_{\Q}$-modules:
\[
\hh = \left(\oplus_{i=1}^r V_i \right) \oplus \left(\oplus_{i=1}^s \Gg_i \right),
\]
where $\oplus V_i = V_{\Q}$ and $\oplus \Gg_i = \Gg_{\Q}$, since both subspaces are
$G_{\Q}$-invariant. For each of these irreducible representations, we can choose
highest weight vectors $v_i \in V_i$ and $w_i \in \Gg_i$, and we then consider
\[
\ul{V}_i := \text{Dist}(\ul{G}_{\Z})v_i, \quad \ul{\Gg}_i := \text{Dist}(\ul{G}_{\Z})w_i,
\]
where $\text{Dist}(\ul{G}_{\Z})$ the algebra of distributions of $\ul{G}_{\Z}$ 
(see \cite[\nopp I.7.7]{Jantzen}). Analogously to
\cite[\nopp II.8.3]{Jantzen}, we have that $V_i = \Q \otimes_{\Z} \ul{V}_i$,
$\Gg_i = \Q \otimes_{\Z} \ul{\Gg}_i$ and that $\ul{V}_{\Z} := \oplus \ul{V}_i$
is a $\ul{G}_{\Z}$-stable lattice inside $V_{\Q}$. As before, set $\ul{V} := \ul{V}_{\Z} \otimes \oO_F$,
which is a $\ul{G}_{\oO_F}$-stable lattice inside $V$.
By scaling the highest weight vectors if necessary, we will assume that $E \in \ul{V}(\oO_F)$.

We can also consider an integral structure $\ul{B}$ on $B$. We can take
the polynomials $p_{d_1},\dots,p_{d_r} \in F[V]^G$ determined in Section \ref{subs:Transverse}
and rescale them using the $\G_m$-action $t \cdot p_{d_i} = t^{d_i}p_{d_i}$
to make them lie in $\oO_F[\ul{V}]^{\ul{G}}$. We let $\ul{B} := \Spec \oO_F[p_{d_1},\dots,p_{d_r}]$
and write $\pi \colon \ul{V} \to \ul{B}$ for the corresponding morphism.
We may additionally assume that the discriminant $\Delta$ defined in
Section \ref{subs:Vinberg} lies in $\oO_F[\ul{V}]^{\ul{G}}$, where the coefficients
of $\Delta$ in $\oO_F[p_{d_1},\dots,p_{d_r}]$ may be assumed to not have a
common divisor.

The following lemma will be convenient to us in the future
(cf. \cite[Lemma 2.8]{ThorneE6}):

\begin{lemma}
There exists a non-zero $N_0 \in \oO_F$ such that for all primes $\pp$ and
for all $b \in \ul{B}(\oO_{\pp})$ we have $N_0 \cdot\kappa_b \in \ul{V}(\oO_{\pp})$.
\end{lemma}

Our arguments in Section \ref{section:construct} will implicitly rely on integral geometric properties of
the representation $(G,V)$. In there, we will need to avoid
finitely many primes, or more precisely to work over $S = \Spec \oO_F[1/N']$
for a suitable $N' \in \oO_F$. By combining the previous lemma and the 
spreading out properties in \cite[\S 7.2]{LagaThesis}, we get:
\begin{prop}
\label{prop:admissible}
There exists a non-zero element $N' \in \oO_F$ such that:
\begin{enumerate}
\item For every $b \in \ul{B}(\oO_F)$, the corresponding Kostant section 
$\kappa_b$ is $G(F)$-conjugate to an element in $\frac{1}{N'}\ul{V}(\oO_F)$.
\item $N'$ is admissible in the sense of \cite[\S 7.2]{LagaThesis}.
\end{enumerate}
\end{prop}
Fix an element $N' \in \oO_F$ satisfying the conclusions of Proposition
\ref{prop:admissible}.
From now on, we will simplify notation by dropping the underline notation
for the objects defined over $\oO_F$, and just refer to $\ul{G},\ul{V},\dots$
as $G,V,\dots$ by abuse of notation.

To end this section, we consider some further integral properties
of the Kostant section. In Section \ref{subs:Transverse}, we considered 
$\kappa$ defined over $F$, and now we will consider some of its properties
over $\oO_{\pp}$. Consider the decomposition
\[
\hh = \bigoplus_{j \in \Z} \hh_j
\]
according to the height of the roots. If $P_H^{-}$ is the negative Borel
subgroup of $H$, $N_H^{-}$ is its unipotent radical and $\mathfrak{p}^-$
and $\mathfrak{n}^-$ are their respective Lie algebras,
we have $\mathfrak{p}^- = \bigoplus_{j \leq 0} \hh_j$, 
$\mathfrak{n}^- = \bigoplus_{j < 0} \hh_j$ and $[E,\hh_j] \subset \hh_{j+1}$. 

\begin{theorem}
\label{theo:intKostant}
Let $R$ be a ring in which $N'$ is invertible. Then:
\begin{enumerate}
\item $[E,\mathfrak{n}_R^-]$ has a complement in
$\mathfrak{p}_R^-$ of rank $\text{rk}_R \, \mathfrak{p}_R^- - \text{rk}_R \, \mathfrak{n}_R^-$;
call it $\Xi$.
\item The action map $N_H^- \times (E + \Xi) \to E + \mathfrak{p}^-$ is
an isomorphism over $R$.
\item Both maps in the composition $E + \Xi \to (E + \mathfrak{p}^-) \git N_H^- \to \hh \git H$
are isomorphisms over $R$.
\end{enumerate}
\end{theorem}
\begin{proof}
See \cite[\S 2.3]{AFV}.
\end{proof}

\begin{obs}
If $R$ is a field of characteristic zero, then $\Xi$ can be taken to be $\fz_{\hh}(f_{0})$ and $E + \Xi$
is the same as the Kostant section considered in Section \ref{subs:Transverse}.
We will abuse notation by referring to both the Kostant section defined
in Section \ref{subs:Transverse} and the section in Theorem \ref{theo:intKostant}
by $\kappa$.
\end{obs}

\section{Constructing orbits}
\label{section:construct}

Given an element $b \in B(\oO_F)$ with discriminant weakly divisible by $I^2$
for a squarefree ideal $I$ avoiding certain bad primes, we will show how to construct
a special integral orbit in $V$ in a way that ``remembers $I$''.
We will start with some technical results.
Assume we have a connected reductive group $L$ over a field $k$, together
with an involution $\xi$. As in Section \ref{subs:Vinberg}, the Lie algebra $\mathfrak{l}$
decomposes as $\mathfrak{l} = \mathfrak{l}(0) \oplus \mathfrak{l}(1)$,
according to the $\pm 1$ eigenspaces of $d\xi$. We also write $L_0$
for the connected component of the fixed group $L^{\xi}$.
\begin{defi}
Let $k$ be algebraically closed. We say a vector $v \in \mathfrak{l}(1)$
is \emph{stable} if the $L_0$-orbit of $v$ is closed and its stabiliser $Z_{L_0}(v)$
is finite. We say $(L_0,\mathfrak{l}(1))$ is stable if it contains stable
vectors. If $k$ is not necessarily algebraically closed, we say $(L_0,\mathfrak{l}(1))$
is stable if $(L_{0,k^s},\mathfrak{l}(1)_{k^s})$ is.
\end{defi}
By \cite[Proposition 1.9]{Thorne16}, the involution $\theta$ defined in Section \ref{subs:Vinberg}
is a stable involution, i.e. $(G,V)$ is stable.

We now prove the analogue of \cite[Lemma 2.3]{RTE8}: the proof is very
similar and is reproduced for convenience.
\begin{lemma}
\label{lemma:etale}
Let $S$ be a $\oO_F[1/N']$-scheme. Let $(L,\xi)$, $(L',\xi')$ be two pairs,
each consisting of a reductive group over $S$ whose geometric fibres are adjoint semisimple of
type $A_1$, together with a stable involution. Then for any $s \in S$ there exists
an étale morphism $S' \to S$ with image containing $s$ and an isomorphism
$L_{S'} \to L_{S'}'$ intertwining $\xi_{S'}$ and $\xi_{S'}'$.
\end{lemma}
\begin{proof}
We are working étale locally on $S$, so we can assume that $L = L'$ and
that they are both split reductive groups. Let $T$ denote the scheme
of elements $l \in L$ such that $\Ad(l) \circ \xi = \xi'$: by \cite[Proposition 2.1.2]{ConradRed},
$T$ is a closed subscheme of $L$ that is smooth over $S$. Since a surjective
smooth morphism has sections étale locally, it is sufficient to show
that $T \to S$ is surjective. Moreover, we can assume that $S = \Spec k$
for an algebraically closed field $k$, since the formation of $T$ commutes
with base change.

Let $A,A' \subset H$ be maximal tori on which $\xi,\xi'$ act as an automorphism
of order 2. By the conjugacy of maximal tori, we can assume that $A = A'$
and that $\xi,\xi'$ define the (unique) element of order 2 in the Weyl group.
Write $\xi = a \xi'$ for some $a \in A(k)$. Writing $a = b^2$ for some
$b \in A(k)$, we have $\xi = b \cdot b \cdot \xi' = b \cdot \xi' \cdot b^{-1}$.
The conclusion is that $\xi$ and $\xi'$ are $H(k)$-conjugate (in fact, $A(k)$-conjugate),
which completes the proof.
\end{proof}
The following lemma is the key technical input in our proof. We remark
the the first part was already implicitly proven in the proof of \cite[Theorem 7.16]{LagaThesis}.
\begin{lemma}
\label{lemma:sl2}
Let $\pp$ be a prime ideal of $\oO_F$ not dividing $N'$, where $N'$ is as
in Proposition \ref{prop:admissible}.
\begin{enumerate}
\item Let $b \in B(\oO_{\pp})$ be an element with $\ord_{\pp} \Delta(b) = 1$,
where $\ord_{\pp} \colon F_{\pp}^{*} \to \Z$ is the usual normalised valuation.
Let $v \in V(\oO_{\pp})$ with $\pi(v) = b$. Then, the reduction mod $\pp$ of $v$
in $V(k_{\pp})$ is regular.
\item Let $b \in B(\oO_{\pp})$ be an element with discriminant weakly divisible
by $\pp^2$. Then,  there exists $g_{\pp} \in G(F_{\pp}) \setminus G(\oO_{\pp})$ 
such that $g_{\pp} \cdot \kappa_b \in V(\oO_{\pp})$.
\end{enumerate}
\end{lemma}
\begin{proof}
Let $v_{k_{\pp}} = x_s+x_n$ be the Jordan decomposition of the reduction of
$v$ in $k_{\pp}$. Then, we have a decomposition $\hh_{k_{\pp}} = \hh_{0,k_{\pp}} \oplus \hh_{1,k_{\pp}}$,
where $\hh_{0,k_{\pp}} = \fz_{\hh}(x_s)$ and $\hh_{1,k_{\pp}} = \text{image}(\Ad(x_s))$.
By Hensel's lemma, this decomposition lifts to $\hh_{\oO_{\pp}} = \hh_{0,\oO_{\pp}} \oplus \hh_{1,\oO_{\pp}}$,
with $\ad(v)$ acting topologically nilpotently in $\hh_{0,\oO_{\pp}}$ and
invertibly in $\hh_{1,\oO_{\pp}}$. As explained in the proof of \cite[Lemma 4.19]{LagaE6},
there is a unique closed subgroup $L \subset H_{\oO_{\pp}}$
which is smooth over $\oO_{\pp}$ with connected fibres and with Lie algebra
$\hh_{0,\oO_{\pp}}$.

For the first part of the lemma, assume that $\oO_{\pp}$ has uniformiser $t$. 
We are free to replace $\oO_{\pp}$ for a
complete discrete valuation ring $R$ with the same uniformiser $t$, containing $\oO_{\pp}$
and with algebraically closed residue field $k$. In this case, the spreading out
properties in \cite[\S 7.2]{LagaThesis} guarantee that the derived group of $L$ is of type
$A_1$. Since the restriction of $\theta$ restricts to a stable involution
in $L$ by \cite[Lemma 2.5]{ThorneThesis}, Lemma \ref{lemma:etale} guarantees that there exists an 
isomorphism $\hh_{0,R}^{der} \cong \fsl_{2,R}$ intertwining the action 
of $\theta$ on $\hh_{0,R}^{der}$ with the action of $\xi = \Ad(\diag(1,-1))$ on $\fsl_{2,R}$.
To show that $v_k$ is regular is equivalent to showing that the
nilpotent part $x_n$ is regular in $\hh_{0,k}^{der}$. The elements $v_k$
and $x_n$ have the same projection in $\hh_{0,k}^{der}$, and given that
$v\in \hh_{0,R}^{der,d\theta = -1}$, its image in $\fsl_{2,R}$ is of the form
\[
\begin{pmatrix}
0 & a \\ b & 0
\end{pmatrix},
\]
with $\ord_R(ab) = 1$ by the spreading out properties in \cite[\S 7.2]{LagaThesis}.
In particular, exactly one of $a,b$ is non-zero when reduced to $k$, and
hence $x_n$ is regular in $\hh_{0,k}^{der}$, as wanted.

For the second part, we return to the case $R = \oO_{\pp}$. There exists
$b' \in B(\oO_{\pp})$ such that $\ord_{\pp}\Delta(b+tb') = 1$, where $t$ 
is a uniformiser of $\oO_{\pp}$. Because the Kostant section is algebraic,
we see that $\kappa_b - \kappa_{b+tb'} \in tV(\oO_{\pp})$, and given that
$\kappa_{b+tb'}$ is regular in $V(k_{\pp})$, we get that $\kappa_b$ is
regular in $V(k_{\pp})$. In particular,
this means that the nilpotent part $x_n$ is a regular nilpotent in $\hh_{0,k_{\pp}}^{der}$.
We now claim that we have an isomorphism $\hh_{0,\oO_{\pp}}^{der} \cong \fsl_{2,\oO_{\pp}}$
intertwining $\theta$ and the previously defined $\xi$ and sending the regular
nilpotent $x_n$ to the matrix
\[
e = 
\begin{pmatrix}
0 & 1 \\ 0 & 0
\end{pmatrix}
\]
of $\fsl_{2,k_{\pp}}$.
Indeed, consider the $\oO_{\pp}$-scheme $X = \Isom((L/Z(L),\theta),(\PGL_2,\xi))$,
consisting of isomorphisms between $L/Z(L)$ and $\PGL_2$ that
intertwine the $\theta$ and $\xi$-actions. Using Lemma \ref{lemma:etale},
we see that étale-locally, $X$ is isomorphic to $\Aut(\PGL_2,\xi)$; in
particular, it is a smooth scheme over $\oO_{\pp}$. By Hensel's lemma \cite[Théorème 18.5.17]{EGA4},
to show that $X$ has a $\oO_{\pp}$-point it is sufficient to show that it has
an $k_{\pp}$-point.

Now, consider the $k_{\pp}$-scheme $Y = \Isom((L/Z(L)_{k_{\pp}},\theta,x_n),(\PGL_2,\xi,e))$
of isomorphisms preserving the $\theta$ and $\xi$-actions which send $x_n$ to $e$:
it is a subscheme of $X_{k_{\pp}}$. Again by Lemma \ref{lemma:etale},
$Y$ is étale locally of the form $\Aut(\PGL_2,\xi,e)$, since $\PGL_2^{\xi}$
acts transitively on the regular nilpotents of $\fsl_2^{d\xi = -1}$ for
any field of characteristic not dividing $N'$. In particular, we see that $Y$ is an 
$\Aut(\PGL_2,\xi,e)$-torsor. In this situation, to see that $Y(k_{\pp})$ is 
non-empty and hence that $X(\oO_{\pp})$ is non-empty, it will suffice to see 
that $\Aut(\PGL_2,\xi,e) = \Spec k_{\pp}$. This follows from the elementary
computation of the stabiliser of $e$ under $\PGL_2^\xi$, which can be seen
to be trivial over any field.

In conclusion, $X(\oO_{\pp})$ is non-empty, meaning that there is an isomorphism
$\hh_{0,\oO_{\pp}}^{der} \cong \fsl_{2,\oO_{\pp}}$ respecting $\theta$ and $\xi$,
and we can make it so that the projection of $v$ in $\fsl_{2,\oO_{\pp}}$ is an
element of the form
\[
\begin{pmatrix}
0 & a \\ bt^2 & 0
\end{pmatrix},
\]
with $a,b \in \oO_{\pp}$ and $a \in 1+t\oO_{\pp}$. Moreover, there exists a morphism
$\varphi \colon \SL_2 \to L^{der}_{F_{\pp}}$ inducing the given isomorphism
$\hh_{0,F_{\pp}}^{der} \cong \fsl_{2,F_{\pp}}$, since $\SL_2$ is simply connected. 
The morphism $\varphi$ necessarily respects the grading, and induces a
map $\SL_2(F_{\pp}) \to L^{der}(F_{\pp})$ on the $F_{\pp}$-points. Consider the 
matrix $g_{\pp}=\varphi(\diag(t,t^{-1}))$: it satisfies the conditions 
of the lemma, and so we are done.
\end{proof}

If $I^2$ divides $\Delta(b)$ weakly, we get for each prime $\pp$ dividing
$I$ an element $g_{\pp}$ as above. Now, consider the adelic element 
$g' \in G(\A_{F,f})$ defined by $g_{\pp}$ at every $\pp$ dividing $I$
and by 1 at every other prime. Recall that by Proposition \ref{prop:clNum},
we can fix elements $\beta_1,\dots,\beta_k$ such that
\[
G(\A_{F,f}) = \coprod_{i=1}^k \lp \prod_{\pp} G(\oO_{\pp}) \rp \beta_i G(F).
\]
Thus, we can write $g' = (h_{\pp})_{\pp}\beta_ig_{I}$ for some $h_{\pp}\in G(\oO_{\pp})$
and $g_{I} \in G(F)$. The element $g_I \cdot \kappa_b$ does not necessarily
lie in $V(\oO_F)$, but rather in
\[
V_{\beta_i} = V(F) \cap \beta_i^{-1} \prod_{\pp} V(\oO_{\pp}).
\]
We see that $V_{\beta_i}$ is a lattice inside $V(F)$, commensurable with
$V(\oO_F)$. In particular, the elements of the different $V_{\beta_1},\dots,V_{\beta_k}$
are all contained in $\frac{1}{N_{bad}} V(\oO_F)$ for some $N_{bad} \in \oO_F$. 
In what follows, we will only consider primes $\pp$ that don't divide
this element $N_{bad}$, and in particular we can assume that $N' \mid N_{bad}$
The choice of $g_I$ is not uniquely determined, it is unique up to the action of
\[
G_{\beta_i} = G(F) \cap \beta_i^{-1} \lp \prod_{\pp} G(\oO_{\pp}) \rp \beta_i.
\]
This is a subgroup of $G(F)$ which is commensurable with $G(\oO_F)$.
We further define the distinguished subspace $W_0 \subset V$ as 
\[
W_0 := \bigoplus_{\substack{a \in \Phi_H/\vartheta \\ \text{ht}(a) \leq 1}} V_a,
\]
where the notation is as in Section \ref{section:ResRoots}. We write an element
$v \in W_0(F)$ as $v = \sum_{\text{ht}(\alpha) = 1} v_{\alpha}X_{\alpha} + \sum_{\text{ht}(\beta) \leq 0} v_{\beta} X_{\beta}$,
where each $X_{\alpha},X_{\beta}$ generates each root space $V_{\alpha},V_{\beta}$
and $v_\alpha,v_\beta \in F$.
Then, we can define the \emph{$Q$-invariant} of $v \in W_0$ as $Q(v) = \prod_{\text{ht}(\alpha) = 1} v_{\alpha}$.
For a squarefree ideal $I \subset \oO_F$, set $G_I = G(F) \setminus \cup_{\pp \mid I}(G(\oO_{\pp}) \cap G(F))$,
that is, the elements of $G(F)$ ``having denominators in $I$''.
Now, define:
\[
W_{i,M} := \left\{v \in V_{\beta_i}\,\middle|\; v = g_I\kappa_b,\, I \text{ squarefree},\, I \text{ coprime to }N_{bad}, \, NI > M,\, g_I \in G_I,\, b \in B(\oO_K)\right\}.
\]
The main result of the section is the following:
\begin{prop}
\label{prop:p-orbit}
Let $b \in B(\oO_F)$, and assume that $\Stab_{G(F)} \kappa_b = \{e\}$.
\begin{enumerate}
\item Let $I$ be a squarefree ideal coprime to $N_{bad}$ and satisfying $NI > M$.
If $I^2$ weakly divides $\Delta(b)$, then $W_{i,M} \cap \pi^{-1}(b)$ is non-empty.
\item If $w \in W_{i,M} \cap W_0$, then $\prod_{v\in M_{\infty}}|Q(w)|_v > M$.
\end{enumerate}
\end{prop}
\begin{proof}
We start by proving the first item. In the above discussion, we constructed
an element $g_I \in G(F)$ such that $g_i \kappa_b \in V_{\beta_i}$.
Given that $I$ is coprime to $N_{bad}$, the construction shows that $g_I \in G_I$,
as otherwise we would have $(\beta_i)_{\pp} \notin G(\oO_{\pp})$ for some prime
$\pp\mid I$, a contradiction.

We now prove the second item. Assume $w = g_I \kappa_b$ for suitable $I$
and $b$. It suffices to show that $\ord_{\pp}Q(w) \geq 1$: in particular,
it suffices to assume that $I = \pp$. Given that $H$ is an adjoint group,
there exists a $t \in T_H(F)$ that makes all the height-one coefficients
of $t \kappa_b$ be equal to one, and in this case we see that actually
$t \in T(F)$. By Theorem \ref{theo:intKostant}, there exists a unique $\gamma \in N_H^{-}(F)$
such that $\gamma t \kappa_b = w$; by taking $\theta$-invariants in the 
isomorphisms of Theorem \ref{theo:intKostant}, we see that 
$\gamma \in N_H^{-,\theta}(F)$. Since the stabiliser is trivial,
we see that $g = \gamma t$, or in other words that 
$g \in P(F) \setminus (P(\oO_{\pp})\cap P(F))$.

Assume, for the sake of contradiction, that $Q(v)$ is invertible in $\oO_{\pp}$,
so that all the height-one coefficients of $w$ are invertible. Then, there
exists a $t' \in T(\oO_{\pp})$ making all the height-one coefficients of
$t'w$ be equal to one, and by Theorem \ref{theo:intKostant}, there exists
at most one element $\gamma'$ in $N^{-}(\oO_{\pp})$ such that $\gamma' t' \kappa_b = w$.
This would imply that $g \in P(\oO_{\pp})$, a contradiction.
\end{proof}

\begin{ex}
Our construction is inspired by the construction in
\cite[Sections 2.2 and 3.2]{BSWsquarefree} for the case $A_n$, for $F = \Q$. In that
case, $C \to B$ corresponds to the family of hyperelliptic curves $y^2 = f(x)$,
where $f(x)$ has degree $n+1$ (there is a slight difference between this
paper and \cite{BSWsquarefree}, in that we consider $f(x)$ without an $x^n$
term while they consider polynomials with a possibly non-zero linear term;
we ignore this difference for now). The main goal of \cite[Sections 2.2 and 3.2]{BSWsquarefree}
is to construct an embedding
\[
\sigma_m \colon \cW_2^{(m)} \to \frac{1}{4}W_0(\Z) \subset \frac{1}{4}V(\Z) ,
\]
where $\sigma_m(f)$ has characteristic polynomial $f$ and $Q(\sigma_m(f)) = m$. 
By taking the usual pinning in $\SL_{n+1}$, we see that our
space $W_0$ corresponds to the space of symmetric matrices in $\fsl_{n+1}$ where the
entries above the superdiagonal are zero, and the height-one entries are
precisely those in the superdiagonal. An explicit section of $B$ can
be taken to lie in $\frac{1}{4}W_0(\Z)$: namely, if $n$ is odd, the matrix
\[
B(b_1,\dots,b_{n+1}) = \begin{pmatrix}
0&1&&&&&&&\\
&0&\ddots&&&&&&\\
&&&1&&&&&\\
&&&0&1&&&&\\
&&&\frac{-b_2}{2}&-b_1&1&&&\\
&&\iddots&-b_3&\frac{-b_2}{2}&0&1&&\\
&\frac{-b_{n-2}}{2}&\iddots&\iddots&&&&\ddots&\\
\frac{-b_{n}}{2}&-b_{n-1}&\frac{-b_{n-2}}{2}&&&&&0&1\\
-b_{n+1}&\frac{-b_{n}}{2}&&&&&&&0
\end{pmatrix}
\]
can be seen to have characteristic polynomial $f(x) = x^{n+1} + b_1x^n +\dots+b_nx +b_{n+1}$.
(if $n$ is even, a similar matrix can be given).
The main observation in this case is that if $m^2$ weakly divides $\Delta(f)$,
then there exists an $l\in \Z$ such that $f(x+l) = x^{n+1} + p_1x^n +\dots+mp_nx +m^2p_{n+1}$
(cf. \cite[Proposition 2.2]{BSWsquarefree}). Then, if $D = \diag(m,1,\dots,1,m^{-1})$,
we observe that the matrix
\[ 
D(B(p_1,\dots,p_{n-1},mp_n,m^2p_{n+1})+lI_{n+1})D^{-1}
\]
is integral, has characteristic polynomial $f(x)$ and the entries in the 
superdiagonal are $(m,1,\dots,1,m)$. Thus, this matrix has $Q$-invariant $m$, as desired.
\end{ex}
\begin{obs}
Our $Q$-invariant is slightly different to the $Q$-invariant defined in
\cite{BSWsquarefree}, which is defined in a slightly more general subspace of $V$. When
restricting to $W_0(F)$, their $Q$-invariant turns out to be a product
of powers of the elements of the superdiagonal, whereas in our case we
simply take the product of these elements. This difference does not affect 
the proof of Theorem \ref{theorem:tail}, and we can also see that for both 
definitions the $Q$-invariant in the previous example is $m$.
\end{obs}

\section{Reduction theory}
\label{section:Reduction}

Before we are able to prove our main results, we need to establish some
results about reduction theory. Mainly, we will be concerned about defining
appropriate heights in the GIT quotient $B$, and constructing a suitable
``box-shaped'' fundamental domain for the action of $\Gamma$ on $G(F_{\infty})$,
where $\Gamma$ will be an arithmetic subgroup of $G(F)$. In the future,
we will use these constructions for $\Gamma = G_{\beta_i}$, where $G_{\beta_i}$
was defined in Section \ref{section:construct}.

\subsection{Heights}
\label{subs:heights}
Recall that $M_{\infty}$ is the set of archimedean places of $F$.
As stated in the introduction, for an element $b \in B(F)$ we define its
\emph{height} to be
\[
\Ht(b) := (NI_b) \prod_{v \in M_{\infty}} \sup\left(|p_{d_1}(b)|_{v}^{1/d_1},\dots,|p_{d_r}(b)|_{v}^{1/d_r}\right),
\]
where $I_b$ is the scaling ideal $\{a \in F \mid a^{d_i}p_{d_i}(b) \in \oO_F, \, \forall i\}$.
We can check that this height is $\G_m(F)$-invariant using the product
formula. A consequence of this is that when $|M_{\infty}| > 1$, the set
of elements of $B(\oO_F)$ having height less than $X$ is infinite. To remedy
that, instead of counting elements in $B(\oO_F)$ we will count the 
number of elements of $\G_m(F) \backslash B(F)$ having height less than 
$X$. We have the following result by Deng (see \cite[Theorem (A)]{Deng}):
\begin{theorem}
\label{theo:heights}
We have
\[
\# \{b \in \G_m(F) \backslash B(F) \mid \Ht(b) \leq X\} = C_B X^{\sum_{i} d_i} + O(X^{\sum_i d_i-\delta_B}),
\]
where $C_B,\delta_B$ are real positive constants.
\end{theorem}
The constants $C_B$ and $\delta_B$ are given explicitly in the statement
of \cite[Theorem (A)]{Deng}.
\begin{obs}
Theorem \ref{theo:heights} is the main reason why we are choosing to work
with this height. There are other natural heights that might be considered,
such as the Weil height:
\[
\Ht_{\text{Weil}}(b) = \prod_{v \in M_F} \sup_i \{|p_{d_i}(b_v)|_v^{1/d_i}\}, 
\]
where now the product is taken over all places of $F$, finite and infinite.
This product is also $\G_m(F)$-invariant by the product formula; however,
we are not aware of any asymptotics for this height in the style of
Theorem \ref{theo:heights}. Having results like that will be very useful
when performing the geometry-of-numbers arguments in Section \ref{section:Counting}.

Moreover, there is a natural interpretation for our choice of height.
In \cite{ESZB}, a natural height on stacks is defined, which in the
particular case of weighted projective spaces turns out to agree with our
choice of height (cf. \cite[\S 3.3]{ESZB}).
\end{obs}
We note the following fact about the quantity $\sum_{i} d_i$ (see \cite[Lemma 8.1]{LagaThesis}):
\begin{prop}
We have $\sum_i d_i = \dim V$.
\end{prop}
For our argument, it will be useful to construct a set $\Sigma \subset B(\oO_F)$ 
which is a fundamental domain for the action of $\G_m(F)$
over $B(F)$, so that it suffices to count elements with invariants in $\Sigma$ to prove
Theorem \ref{theorem:tail}. We will do so following \cite[\S 3.4]{BSWcoregular}.
For a finite prime $\pp$, we set
\[
B'(\oO_{\pp}) = \{b \in B(\oO_{\pp}) \mid v_{\pp}(p_{d_i}(b)) < d_i \text{ for some }i\}.
\]
For every $b_{\pp} \in B(F_{\pp})$, there exists $g_{\pp} \in \G_m(F_{\pp})$
such that $g_{\pp} b_{\pp} \in B'(\oO_{\pp})$; and moreover this
$g_{\pp}$ is unique up to the action of $\G_m(\oO_{\pp})$. In a similar
spirit to Section \ref{section:construct}, for any $\gamma \in \cl(\G_m)$
we can consider the sets
\begin{align*}
B_{\gamma} &= B(F) \cap \gamma^{-1} \prod_{\pp} B'(\oO_{\pp}); \\
\G_{m,\gamma} &= \G_m(F) \cap \gamma^{-1} \prod_{\pp} \G_m(\oO_{\pp}) \gamma = \oO_F^{\times}.
\end{align*}
Let $\gamma_{1},\dots,\gamma_{k}$ be representatives of $\cl(\G_m)$, 
which is finite by Proposition \ref{prop:clNum} (in fact, it coincides
with the class group of $F$). Then, we have a
bijection
\[
\coprod_{i=1}^k \oO_F^{\times}\backslash B_{\gamma_i} \longleftrightarrow \G_m(F) \backslash B(F),
\]
given by inclusion. Indeed, if for some $g \in \G_m(F)$ we have
$v_1 \in B_{\gamma}$ and $g v \in B_{\gamma'}$, then for all primes $\pp$
we have that $\gamma_{\pp} v \in B'(\oO_{\pp})$ and $\gamma_{\pp}' g v \in B'(\oO_{\pp})$.
This implies that $\gamma_{\pp}'g\gamma_{\pp}^{-1} \in \G_m(\oO_{\pp})$
which means that $\gamma = \gamma'$, from which injectiveness follows. The map
is surjective by construction.

We can modify our choices of representatives $\gamma_i$ of $\cl(\G_m)$
so that $B_{\gamma_i} \subset B(\oO_F)$, simply by multiplying by appropriate
elements of $\G_m(F)$. Moreover, for a given choice of $\gamma_i$, the
ideal $I_b$ is independent of the choice of $b \in B_{\gamma_i}$: in
fact, a computation shows that $I_b$ is equal to the ideal corresponding
to $\gamma_i$, regardless of the choice of $b$ (as long as $b\neq 0$).

To construct the fundamental domain $\Sigma$, it suffices to construct
fundamental domains $\Sigma_i$ in each
$B_{\gamma_i}$ separately. Given that in the future we will want to impose
congruence conditions in $\Sigma_i$, we will also define $\Sigma_i$ itself 
as a set defined by congruence conditions: that is, defined as the 
intersection of sets $\Sigma_{i,\pp} \subset B(\oO_{\pp})$ for all primes 
$\pp$ (finite and infinite).

For finite primes $\pp$, the subset $\Sigma_{i,\pp}$ is given by 
$\gamma_{\pp}^{-1} B'(\oO_{\pp})$. For infinite primes, we consider the set $F_{\infty} = \prod_{v \in F_{\infty}} F_v$ and
the subset $F_{\infty}^1$ of $F_{\infty}$ consisting of those elements
$(\alpha_v)_{v \in M_{\infty}}$ such that $\prod_v |\alpha_v|_v = 1$. We
further consider $\Lambda$, a compact subset of $F_{\infty}^1$ such that
$F_{\infty}^1 = \Lambda \oO_F^{\times}$.  Let $B(1)$ be the set of elements
$b \in B(F_{\infty})$ such that, for all $v \in M_{\infty}$:
\[
\sup_i\{|p_{d_i}(b)|_{v}^{1/d_i}\} = \Ht(b)^{1/|M_{\infty}|}.
\]
Let $\overline{B}(1)$ be a measurable set with boundary of measure 0 which
is a fundamental domain for the action of the roots of unity $\mu_F$ over
$B(1)$. Then, the set $\Lambda\cdot(\overline{B}(1))$ is a fundamental domain
for the action on $\oO^{\times}$ over $B(F_{\infty})$, which we take as
our $\Sigma_{i,\infty}$. In conclusion, we obtain our fundamental domain
$\Sigma_i$ by combining the congruence conditions $\Sigma_{i,\pp}$ and $\Sigma_{i,\infty}$.

The fact that $\Sigma \cap B(F)_{X}$ is finite follows from the fact
that $\Sigma_{\infty}$ consists of elements whose local heights differ
at most by an absolute constant. Thus, if an element has bounded height, then each
of the local heights has to be bounded, yielding a finite number of elements
overall.

\subsection{Measures on $G$}
\label{subs:measures}

Recall that $\Phi_G = \Phi(G,T)$ is the set of roots of $G$. The Borel
subgroup $P^+$ of $G$ determines a root basis $S_G$ and a set of
positive/negative roots $\Phi_G^{\pm}$, compatible with the choice of positive
roots in $H$ determined by the pinning of Section \ref{subs:Vinberg}. 
Let $N$ be the unipotent radical of the negative Borel subgroup 
$P$. We will make an appropriate choice of maximal compact
subgroup of $G(F_{\infty})$:
\begin{lemma}
\label{lemma:cartan}
There exists a maximal compact subgroup $K$ of $G(F_{\infty})$ such that
$T(F_{\infty}) \cap K = \{t = (t_1,\dots,t_r) \in T(F_{\infty}) \mid |t_i|_v = 1, \, \forall v \in M_{\infty}, \forall i\}$.
\end{lemma}
\begin{proof}
We will choose a maximal compact $K_v \subset G(F_v)$ for every infinite
place $v$ satisfying that $T(F_v) \cap K_v = \{t = (t_1,\dots,t_r) \in T(F_{v}) \mid |t_i|_v = 1, \, \forall i \}$.
Then, it will be enough to define $K = \prod_v K_v$.

Assume that $F_v = \R$. By the Isomorphism Theorem (see e.g. \cite[Theorem 6.1.17]{ConradRed}),
the involution corresponding to $-1$ in the root system $\Phi_G$ induces
an involution $\vartheta \colon G(F_v) \to G(F_v)$ that acts as inversion
in $T(F_v)$. Moreover, this involution $\vartheta$ is a Cartan involution:
to check this, we need to verify that the form $B_{\vartheta}(X,Y) := -B(X,\vartheta(Y))$
is positive definite (here, $B$ is the Killing form). If we fix basis elements
$X_{\alpha} \in \Gg_{\alpha}$ for each $\alpha \neq 0$ in $\Phi_G$, we can
additionally require that $d\vartheta(X_{\alpha}) = -X_{-\alpha}$. Then,
it is straightforward to check that $B_{\vartheta}$ is positive definite.

Now, assume that $F_v = \C$. An involution $\vartheta \colon G \to G$ is a Cartan
involution if and only if $\{g \in G(\C) \mid \vartheta(g) = \ol{g}\}$
is a maximal compact subgroup of $G(\C)$. Note that in the split torus
$T$, the involution $\vartheta_T \colon T \to T$ sending $t \in T(\C)$
to $t^{-1}$ is a Cartan involution: if $t = (t_1,\dots,t_k) \in (\C^{\times})^k$,
then $\vartheta_T(t) = \ol{t}$ if and only if $|t_i| = 1$ for all $i$, so
the set $\{t \in T(\C) \mid \vartheta_T(t) = \ol{t}\}$ is a maximal compact
subgroup of $T(\C)$.
By \cite[Theorem 3.13(1)(b)]{AT}, the Cartan involution $\vartheta_T$
extends to a Cartan involution $\vartheta \colon G \to G$. By taking
$K_v = \{g \in G(\C) \mid \vartheta(g) = \ol{g}\}$, we see that $T(F_v) \cap K_v$
corresponds exactly to those elements with modulus 1, as wanted.
\end{proof}

Consider the subgroup $A \subset T(F_{\infty})$ of elements $a = (a_v)_{v \in M_{\infty}}$
such that $a_v \in \R_{> 0}$ for all $v \in M_{\infty}$. Then, the map
\[
N(F_{\infty}) \times A \times K \to G(F_{\infty})
\]
given by $(n,t,k) \mapsto ntk$ is a diffeomorphism.
The following result is a well-known property of the Iwasawa decomposition:
\begin{lemma}
\label{lemma:IwasawaInt}
Let $dn,dt,dk$ be Haar measures on $N(F_{\infty}),A,K$,
respectively. Then, the assignment
\[
f \mapsto \int_{n \in N(F_{\infty})}\int_{t \in A}\int_{k \in K} f(ntk) |\delta(t)|^{-1} dn\, dt\, dk
\]
defines a Haar measure on $G(F_{\infty})$. Here, 
$\delta(t) = \prod_{\beta \in \Phi_G^{-}} \beta(t) = \det \Ad(t)|_{\Lie \overline{N}(F_{\infty})}$
is an algebraic character obtained from the action of $T$ on $N$, and
$|\delta(t)| = \prod_{v \in M_{\infty}}|\delta(t_v)|_v$.
\end{lemma}

We get the Haar measure $dt$
from the isomorphism with $\prod_{v \in M_{\infty}}(\R_{> 0})^{\# S_G}$,
where $\R_{> 0}$ is given the standard Haar measure $d^{\times}\lambda = d\lambda/\lambda$.

\subsection{Fundamental sets}
\label{subs:FDomains}

In this section, given an arithmetic subgroup $\Gamma \subset G(F)$, we
will construct an exact fundamental domain $\cF$ for the action of $\Gamma$
on $G(F_{\infty})$.
\begin{defi}
A \emph{Siegel set} is a set of the form $\cS = \cup_{i=1}^n \alpha_i\omega_i A_c K$, where
$\alpha_i \in G(F)$, the sets $\omega_i \subset N(F_{\infty})$ are compact and
\[
T_c = \{t \in T(F_{\infty}) \mid |\alpha(t)| \leq c, \; \forall \alpha \in S_G\}, \quad A_c = T_c \cap A.
\]
\end{defi}
\begin{obs}
In fact, we will consider subsets of the form $\cS = \cup_{i=1}^n \alpha_i\omega_i A_c' K_1$,
where $A_c'$ and $K_1$ are appropriate subsets of $A_c$ and $K$, respectively.
We will still call such a set a Siegel set, for simplicity.
\end{obs}
We will require our fundamental domain $\cF$ to be ``box-shaped at infinity'',
in the following sense:
\begin{defi}
\label{defi:box}
We say a fundamental domain $\cF$ for the action of $\Gamma$ on $G(F_{\infty})$
is \emph{box-shaped at infinity} if there exist Siegel sets $\cS_1 \subset \cF \subset \cS_2$
such that
\begin{itemize}
\item $\cS_1$ and $\cS_2$ have the same cusps, i.e. both sets are defined
as $\cS_1 = \cup_{i=1}^n \alpha_i \omega_i A_{c_1}' K_1$ and $\cS_2 = \cup_{i=1}^n \alpha_i \omega_i A_{c_2}' K_1$
for the same choice of elements $\alpha_i \in G(F)$ in each case, and the
same choice of subsets $A_c'$ and $K_1$.
\item There exists an open subset $\cU_1 \subset \cS_1$ of full measure
such that each $\Gamma$-orbit on $G(F_{\infty})$ intersects $\cU_1$ at
most once.
\item Every $\Gamma$-orbit on $G(F_{\infty})$ intersects $\cS_2$ at least
once.
\item For sufficiently small $c$, we have that 
$\cS_1 \cap \lp \cup_{i=1}^n \alpha_i  N(F_{\infty}) A_{c} K \rp = \cS_2 \cap \lp \cup_{i=1}^n \alpha_i N(F_{\infty}) A_{c} K \rp$.
\end{itemize}
\end{defi}
In what follows, it will be sufficient to construct $\cS_1$ and $\cS_2$
to obtain $\cF$ due to the following lemma (cf. \cite[Lemma 7]{SSSVcusp}):
\begin{lemma}
\label{lemma:box}
Let $B(G)$ be the Borel $\sigma$-algebra of $G(F_{\infty})$. Assume we
have $\cS_1$ and $\cS_2$ in $B(G)$ such that the maps $\cS_1 \to \Gamma \backslash G(F_{\infty})$
and $\cS_2 \to \Gamma \backslash G(F_{\infty})$ are injective and surjective,
respectively. Then, there is a fundamental domain $\cF$ in $B(G)$ for
the action of $\Gamma$ in $G(F_{\infty})$ such that $\cS_1 \subset \cF \subset \cS_2$.
\end{lemma}
\begin{proof}
Since $\Gamma$ is a discrete subgroup of $G(F_{\infty})$, we can find a
non-empty open set $U \subset G(F_{\infty})$ such that $U^{-1}U \cap \Gamma = \{\id\}$.
Given that $G(F_{\infty})$ is second countable, we can find countably
many elements $\{g_n\}_{n\in \N} \subset G(F_{\infty})$ such that
$G(F_{\infty}) = \cup_{n=1}^{\infty}g_n U$. Let $\cS_3 = \cS_2 \setminus \Gamma \cS_1$,
and define:
\[
\cF_0 = \bigcup_{n=1}^{\infty} \lp g_nU \cap \cS_3 \setminus \bigcup_{i < n} \Gamma(g_i U \cap \cS_3)  \rp.
\]
Finally, set $\cF = \cS_1 \cup \cF_0$, a disjoint union. Then, $\cF \in B(G)$
and the map $\cF \to \Gamma\backslash G(F_{\infty})$ can easily be seen
to be bijective, as wanted.
\end{proof}

\subsubsection{Constructing $\cS_1$}
\label{subs:corners}

To obtain the domain $\cS_1$, we will use general properties of the
Borel-Serre compactification following \cite{BScorners}.

First, consider the Weil restriction of scalars $G' = \Res_{F/\Q} G$,
which is a semisimple group defined over $\Q$. We have an isomorphism
$\phi_G \colon G'(\Q) \cong G(F)$, inducing $\phi_G \colon G'(\R) \cong G(F_{\infty})$. Even though the
group $G$ is split with maximal torus $T$, the group $G'$ will
\emph{never} be split (unless $F = \Q$). Let us denote $T' = \Res_{F/\Q} T$, a maximal
torus of $G'$ and $T_{split}' \subset T'$ the maximal $\Q$-split torus
inside $T'$.

Recall that the $F$-split torus $T$ is obtained by base-changing a $\Q$-split
torus $T_{\Q}$ to $F$. In particular, there is an isomorphism over $\Q$
between the split tori $T_{split}' \cong T_{\Q}$, inducing an isomorphism
of character groups $X^{*}(T_{split}') \cong X^{*}(T_{\Q})$. By \cite[(6.21)]{BT}
this isomorphism induces a bijection between the root systems $\Phi(G,T)$
and $\Phi(G',T_{split}')$. Denote by $S_{G'}$ a positive root basis
of $\Phi(G',T_{split}')$, chose compatibly with the positive root
basis of $\Phi(G,T)$. Let us define
\[
T_c' = \{t \in T_{split}'(\R) \mid |\beta(t)| \leq c,\, \forall \beta \in S_{G'} \}.
\]
Let $P' = \Res_{F/\Q} P$ be a minimal parabolic subgroup of $G'$ containing
$T'$, and define
\[
{}^0P'(\R) = \bigcap_{\chi \in X^{*}(P')} \ker \chi^2.
\]
We then have that $P'(\R) = {}^0P'(\R) \rtimes A_{split}$ by \cite[Proposition 1.2]{BScorners},
where $A_{split}$ denotes the connected component of $T_{split}'(\R)$.
\begin{prop}
\label{prop:compare}
Let $c > 0$. Then, $\phi_G(T_c) \subset {}^0P'(\R) \rtimes (T_{c'}' \cap A_{split})$
for some positive constant $c'$. Similarly, $\phi_G^{-1}(T_c') \subset T_{c'}$
for some $c'$.
\end{prop}
\begin{proof}
Let us begin by recalling the isomorphism $T_{\Q} \cong T_{split}'$. When
base-changing to $\R$, we see that $\phi_G^{-1}(T_{split}'(\R))$ corresponds
to the $t = (t_1,\dots,t_k) \in T(F_{\infty})$ such that each $t_i = (t_{i,v})_{v\in M_{\infty}} \in F_{\infty}$
satisfies $t_{i,v} = t_{i,v'} \in \R^{\times}$ for all infinite places $v,v'$.

Let $t \in T_c$. We can decompose $t = t_1t_2$ with $t_1 \in \phi_{G}^{-1}(T_{split}'(\R))$,
scaled in such a way that $\prod_v |\alpha(t_v)|_v = \prod_v |\alpha(t_{1,v})|_v$
for all $\alpha \in X^{*}(T)$. By \cite[(6.20)]{BT}, there is an isomorphism
$X^{*}(T)_F \xrightarrow{F} X^{*}(T')_{\Q}$ such that for all $\chi \in X^{*}(T)$
and $g \in G(F_{\infty})$, we have that $\prod_v |\chi(g_v)|_v = |(F \circ \chi)(\phi_G(g))|$.
In particular, we have that $\phi(t_2) \in {}^0P'(\R)$. Finally, for
any $t_2 \in \phi_{G}^{-1}(T_{split}'(\R))$, we have that $t_2 \in T_c$
if and only if $\phi_G(t_2) \in T_{c^{1/[F:\Q]}}'$, so choosing
$c' = c^{1/[F:\Q]}$ concludes the proof of the first inclusion. The second
inclusion is analogous.
\end{proof}

We will denote by $K'$ the restriction of the maximal compact $K$
inside $G'$. Now, consider the symmetric space $X = G'(\Q)/K'$. 
For each parabolic $\Q$-subgroup $Q$ of $G'$, let
$S_Q := (R_dQ/(R_uQ \cdot R_dG'))$, where $R_u$ denotes the unipotent radical
and $R_d$ denotes the $\Q$-split part. Then, $S_Q$ is a $\Q$-split torus,
and we let $A_Q := S_Q(\R)^{\circ}$. There is a natural action of $A_Q$ on $X$;
called the geodesic action (see \cite[(3.2)]{BScorners}). Set $e(Q) = A_Q \backslash X$,
and consider
\[
\overline{X} = \coprod_{P \subset G'\text{ parabolic}} e(Q),
\]
which by \cite[(7.1)]{BScorners} naturally has a structure of a manifold
with corners. The topology of $\overline{X}$ is studied in \cite[\S 5, \S 6]{BScorners};
in particular, it is shown that for any parabolic group $Q$, the subset
$X(Q) = \coprod_{R \supset Q} e(R)$ is an open subset of $\overline{X}$.
Taking $Q = G$, we see that $e(G) = X$ is an open submanifold of $\overline{X}$.

Now, let us return to considering the minimal parabolic 
$P' = \Res_{K/\Q}P$, and consider
\[
U_{x,P',c} = {}^0 P'(\R) (A_{P',c} \cdot x),
\]
where $x \in X$ and $A_{P',c} = A_{split} \cap T_c'$, as defined in the
beginning of the section.
The closure $\overline{U_{x,P',c}}$ in $\overline{X}$ is a neighbourhood
of the closure of $e(P')$ in $\overline{X}$. The arithmetic subgroup
$\Gamma$ of $G(F)$ restricts to an arithmetic subgroup $\Gamma'$ of $G'(\Q)$:
when we consider the action of $\Gamma'$ in these sets $U_{x,P',c}$, it
is useful to consider the following result (see \cite[Proposition 10.3]{BScorners}):
\begin{prop}
\label{prop:corners}
There exists $c > 0$ satisfying that for any $g_1, g_2 \in \overline{U_{x,P',c}}$, if there
exists $\gamma \in \Gamma'$ such that $g_1 = \gamma g_2$, then $\gamma \in P' \cap \Gamma'$.
\end{prop}
For our construction of $\cS_1$, we will have to worry about different
cusps at once. Let $\{\alpha_1,\dots,\alpha_m\}$ be a set of representatives
for the double cosets of $\Gamma \backslash G(F) / P(F)$,
which is equivalent to choosing a set of representatives $\{\alpha_1',\dots,\alpha_m'\}$
for the double cosets of $\Gamma' \backslash G'(\Q) / P'(\Q)$.
The interaction of different cusps can be controlled as follows (see \cite[Proposition 10.4]{BScorners}):
\begin{prop}
\label{prop:cusps}
Let $Q,R$ be parabolic subgroups of $G'$, let $x,y \in X$ and let $g \in G'(\Q)$.
Then, $g\cdot U_{x,Q,c} \cap U_{y,R,c} \neq \emptyset$ for all $c > 0$
if and only if $gQg^{-1} \cap R$ is parabolic.
\end{prop}
In particular, when $Q$ is a minimal parabolic subgroup and $Q = R$,
we have that $gQg^{-1} \cap Q$ is parabolic exactly when $g \in Q$.

We will take our set $\cS_1$ in $G(F_{\infty})$ to be contained in a set
of the form $\cup_{i=1}^m \alpha_i \omega A_c K$, where $\alpha_i \in G(F)$
are as above. We will choose $c$ to be ``big enough'':
\begin{prop}
There exists a small enough $c > 0$ such that if
$g_1 \in \alpha_i N(F_{\infty})A_c K$ and $g_2 \in \alpha_j N(F_{\infty})A_c K$
are $\Gamma$-equivalent, then $\alpha_i = \alpha_j$.
\end{prop}
\begin{proof}
By Proposition \ref{prop:compare}, we have that $\phi_G(N(F_{\infty})A_c K) \subset {}^{0}P'(\R)T_{c'}'K$
for some $c'>0$. We choose $c > 0$ small enough so that $c'$ satisfies
the conclusions of Proposition \ref{prop:corners}. Then, if there
exists $\gamma \in \Gamma$ such that $\gamma g_1 = g_2$, then we would
have that  $\alpha_j^{-1} \gamma \alpha_i \in P(F)$,
or in other words that $\alpha_i \in \Gamma \alpha_j P(F)$, which can
only happen if $\alpha_i = \alpha_j$.
\end{proof}
Therefore, we can choose $c$ small enough so that there are no intersections
between the cusps. Now, we need to choose appropriate subsets of $N(F_{\infty})$,
$A_c$ and $K$ so that each cusp contains at most one $\Gamma$-orbit.
\begin{lemma}
\label{lemma:omega}
There exist compact sets $\omega_i$ inside $N(F_{\infty})$ such that
every $\alpha_i^{-1} \Gamma \alpha_i \cap N$-orbit intersects $\omega_i$,
and also satisfying that inside a set of full measure, every $\alpha_i^{-1} \Gamma \alpha_i \cap N$-orbit
intersects $\omega_i$ exactly once.
\end{lemma}
\begin{proof}
First, we recall that for any non-zero root $\beta \in \Phi_G$, there
is an isomorphism $u_{\beta}\colon \G_a \to U_{\beta}$, for some subgroup $U_{\beta} \subset G$.
Let $\beta_1,\dots,\beta_m$ be an ordering of the roots of $\Phi_G^-$,
such that $|\Ht(\beta_i)| \leq |\Ht(\beta_{i+1})|$. Then, by 
\cite[Theorem 5.1.13]{ConradRed}, there is an isomorphism of algebraic
varieties (not necessarily a group morphism)
\[
\prod_{i=1}^m U_{\beta_i} \to N
\]
defined over $\Z$, which is just given by the product map. Given that 
$N_{i} := \alpha_i^{-1} \Gamma \alpha_i \cap N$ is an arithmetic
subgroup of $N(F)$, we can arrange everything so that the elements of
$N_{i}$ correspond precisely to those elements $u_{\beta_1}(x_1)\dots u_{\beta_m}(x_m)$
such that $x_m \in \oO_F$. Choose a compact subset $\Lambda$ of $F_{\infty}$ such
that $\oO_F + \Lambda = F_{\infty}$, such that the interior of $\Lambda$
has full measure inside $\Lambda$, and such that no two elements in the
interior differ by an element of $\oO_F$. Now, set
\[
\omega_i = \{u_{\beta_1}(x_1)\cdots u_{\beta_m}(x_m) \mid x_1,\dots,x_m \in \Lambda\}.
\]
Using the commutator relations in \cite[Proposition 5.1.14]{ConradRed},
it is not difficult to see that $\omega_i$ satisfies the conclusions of
the lemma.
\end{proof}
Now, we worry about the action of $T$ on $A_c$ and $K$. More specifically,
we need to account for the action of the group $T_i := \alpha_i^{-1}\Gamma\alpha_i \cap T(F)$,
which is an arithmetic subgroup of $T(F)$. In particular, it has to be
commensurable to $T(\oO_F)$, which means in particular that $T_i \subset T(\oO_F)$.

Let $w_F \subset \oO_F^*$ be the subgroup of roots of unity of $\oO_F$.
Under the identification $T(\oO_F) \cong (\oO_F^{*})^{\# S_G}$,
consider the subgroup $T_{w,i} \subset T_i$ corresponding to the appropriate
subgroup of $(w_F)^{\# S_G}$ lying in $T_i$. Then, we know that $T_{w,i} \subset K$ by Lemma \ref{lemma:cartan}.
We consider a fundamental domain for the left action of $T_{w,i}$ on $K$,
which we will denote $K_1$.

Let $|T_i|$ denote the subset of $A(F_{\infty})$ which is the image of
$T_i$ under the projection map that sends $t = (t_1,\dots,t_k) \in T_i$
to $(a_1,\dots,a_k)$, where $a_{j,v} = |t_j|_v$. We denote a fundamental
domain for the action of $|T_i|$ on $A_c$ by $A_c'$.

Finally, we let $\cS_1 = \cup_{i=1}^m \alpha_i \omega_i A_c' K_1$.
\begin{theorem}
There exists an open subset $\cU_1$ of $\cS_1$ of full measure such that
any $\Gamma$-orbit in $G(F_{\infty})$ intersects $\cU_1$ at most once.
\end{theorem}
\begin{proof}
Let $\omega_i' \subset \omega_i$ be the subset of full measure described in Lemma \ref{lemma:omega},
and let
\[
\cU_1 = \cup_{i=1}^m \alpha_i \omega_i' A_c' K_1.
\]
Now, let $g_1,g_2 \in \cU_1$, and let $\gamma \in \Gamma$ be an element
such that $\gamma \cdot g_1 = g_2$. By Proposition \ref{prop:cusps}, we know that $g_1$ and $g_2$
have to lie in the same cusp; that is, we have $g_1,g_2 \in \alpha_i \omega_i' A_c' K_1$.
Write $g_1 = \alpha_i n_1t_1k_1$ and $g_2 = \alpha_i n_2 t_2 k_2$.
We have that by Proposition \ref{prop:corners} that $\alpha_i^{-1}\gamma \alpha_i \in P$,
or in other words that $\gamma \in \alpha_i P \alpha_i^{-1} \cap \Gamma$.
Let $\gamma = \alpha_i n_0t_0 \alpha_i^{-1}$, where $n_0\in \alpha_i^{-1}\Gamma\alpha_i \cap N(F)$
and $t_0 \in \alpha_i^{-1}\Gamma\alpha_i \cap T(F)$. Then, the condition that $\gamma \cdot g_1 = g_2$ becomes
\[
\alpha_i n_0 (t_0 n_1 t_0^{-1}) t_0 t_1 k_1 = \alpha_i n_2 t_2 k_2.
\]
Let us write $t_0 = t_at_k$, where $t_a \in A(F_{\infty})$ and $t_k \in K$.
Then, the uniqueness in the Iwasawa decomposition gives us $n_0 t_0 n_1 t_0^{-1} = n_2$,
$t_at_1 = t_2$ and $t_k k_1 = k_2$. By construction of $A_c'$, we must
have that $t_a = 1$, and therefore that $t_k \in T_w$, so by construction
of $K_1$ we also have that $t_k = 1$. Then $t_0 = 1$, so the equation
$n_0 n_1 = n_2$ also gives $n_0 = 1$ by construction of $\omega_i'$.
\end{proof}

\subsubsection{Constructing $\cS_2$}

Having chosen $\cS_1$, we will now construct a compatible $\cS_2$.

\begin{prop}
\label{prop:cS2}
There exists $c > 0$ such that $G(F_{\infty}) = \Gamma \cdot \cup_{i=1}^m \alpha_i \omega_i A_c' K_1$,
where $\alpha_i$, $\omega_i$, $A_c'$ and $K_1$ are as in Section \ref{subs:corners}.
\end{prop}
\begin{proof}
Using Proposition \ref{prop:compare} and restriction of scalars, the results
in \cite{SpringerGF} tell us that there exists a $c > 0$ such that
\[
G(F_{\infty}) = \Gamma \cdot \bigcup_{i=1}^m \alpha_i \omega' T_c K
\]
for some compact subset $\omega' \in N(F_{\infty})$. Now, assume that
we have $g \in G(F_{\infty})$ written as $g = \alpha_i n_0 t_0 k_0$ for
some $n_0 \in \omega'$, $t_0 \in T_c$ and $k_0 \in K$. We want to see
that $g \in \Gamma \cdot \cup_{i=1}^m \alpha_i \omega_i A_c' K_1$.

By construction, we will have that $t_0 k_0 = t_{\Gamma} t_1 k_1$, where
$t_{\Gamma} \in \alpha_{i}^{-1}\Gamma \alpha_i \cap T(F)$, $t_1 \in A_c'$ and $k_1 \in K_1$. We also have that there
exists $n_{\Gamma} \in \alpha_i^{-1} \Gamma \alpha_i \cap N(F)$ and $n_1 \in \omega_i'$
such that $n_{\Gamma}t_{\Gamma} n_0 t_{\Gamma}^{-1} = n_1$. Then, we have that
\[
(\alpha_i n_{\gamma} t_{\Gamma} \alpha_i^{-1})\cdot \alpha_i n_0 t_0 k_0 = \alpha_i n_{\gamma} t_{\Gamma} n_0 t_{\Gamma}^{-1} t_{\Gamma} t_1 k_1 = \alpha_i n_1 t_1 k_1,
\]
as wanted.
\end{proof}

We choose $\cS_2$ to be the set constructed in Proposition \ref{prop:cS2}.
It is clear that $\cS_1$ and $\cS_2$ satisfy the conditions stated in
Definition \ref{defi:box}. Therefore, by Lemma \ref{lemma:box}, 
we have constructed the desired box-shaped fundamental domain $\cF$.

For the future, it will be useful to record the following property of $\cS_2$
(known as the Siegel property):
\begin{prop}
\label{prop:Siegel}
The size of the fibres of the map $\cS_2 \to \Gamma \backslash G(F_{\infty})$
is uniformly bounded.
\end{prop}
\begin{proof}
It suffices to show that the set $\{\gamma \in \Gamma \mid \gamma \cS_2 \cap \cS_2 \neq \emptyset\}$
is finite. For $F = \Q$, this follows from \cite[Corollaire 15.3]{Bor69},
and in the general case we reduce to $F = \Q$ using restriction of scalars.
\end{proof}

\section{Counting reducible orbits}
\label{section:Counting}

Before we are able to prove Theorem \ref{theorem:tail}, we need to obtain
an estimate on the number of reducible $\Gamma$-orbits on some $\Gamma$-invariant
lattice $V_0$ of $V(F_{\infty})$ which is commensurable with $V(\oO_F)$.
Given that, as observed in Section \ref{subs:heights}, there might be
infinitely many orbits with bounded height, we will restrict to the elements
with invariants lying on the fundamental domain $\Sigma \subset B(\oO_K)$.
Recall that $\Sigma$ was the disjoint union of finitely many sets $\Sigma_i$,
where for each $b \in \Sigma_i$ we had
\begin{equation}
\label{eq:height}
\Ht_i(b) = C_i\prod_{v \in M_{\infty}} \sup\{ |p_{d_i}(b)|_v^{1/d_i}\},
\end{equation}
for some constant $C_i$ only dependent on $i$ (not on the choice of $b$
inside $\Sigma_i$). We will count the number of reducible $\Gamma$-orbits
in $V_0$ having invariants in $\Sigma_i$ for each $i$ separately. For
simplicity, we will denote $\Sigma_i$ simply by $\Sigma$ in this section.

For any element $b \in B(F_{\infty})$, we \emph{define} its height to be
given by the expression in \eqref{eq:height}, and for any $v \in V(F_{\infty})$
we also define $\Ht(v) := \Ht_i(\pi(v))$.

Let $\Lambda$ be the embedding of $\R_{> 0}$ inside $F_{\infty}$ given by
sending $x \in \R_{> 0}$ to $|M_{\infty}|$ copies of $x$ inside every 
infinite place of $F$. For an element $\lambda \in \Lambda$ and $v \in V(F_{\infty})$, we have that
$\Ht(\lambda v) = \lambda^{[F:\Q]}\Ht(v)$, or in other words that the
function $\Ht$ is homogeneous of weight $[F:\Q]$.

We will prove the following:
\begin{theorem}
\label{theo:constant}
Let $N(\Gamma,V_0,X)$ denote the number of reducible $\Gamma$-orbits in
$V_0$ having height less than $X$ and invariants lying in $\Lambda\Sigma$. Then,
we have
\[
N(\Gamma,V_0,X) = CX^{\dim V} + O(X^{\dim V - \delta}).
\]
The constant $C$ depends only on $\Gamma$ and $V_0$, and the constant
$\delta$ can be chosen independently of $\Gamma$ and $V_0$.
\end{theorem}

\subsection{Averaging and reductions}
\label{subs:reductions}

By analogous arguments to \cite[\S 2.9]{ThorneE6}, there exist open subsets $L_1,\dots,L_k$
covering $\{b \in B(F_{\infty}) \mid \Ht(b) = 1,\, \Delta(b) \neq 0\}$ such that
for a fixed $i$, the quantity $r_i := \# \Stab_G(v)(F_{\infty})$ remains constant for 
any choice of $v \in \pi^{-1}(L_i)$. Let us denote by $\Lambda$ the embedding
of $\R_{> 0}$ inside $F_{\infty}$ that sends $x \in \R_{> 0}$ to $(x,\dots,x) \in F_{\infty}$,
and denote $V_i := V_0^{red} \cap G(F_{\infty}) \Lambda\kappa(L_i)$.
Fix a compact left and right $K$-invariant set $G_0 \subset G(\R)$ which is the
closure of a non-empty open set, for which we assume that $G_0 = G_0^{-1}$.
An averaging argument just as in \cite[\S 2.3]{BSquartics}) yields
\begin{equation}
\label{eq:averaging}
N(\Gamma,V_i,X) = \frac{1}{r_i \vol(G_0)} \int_{g \in \cF} \# \{v \in V_0^{red} \cap (g G_0 \Lambda \kappa(L_i \cap \Sigma))_{< X} \} dg. 
\end{equation}
Here, the subscript $< X$ means we are restricting to elements of height
less than $X$. For simplicity, we will denote $\cB_X = (G_0 \Lambda \kappa(L_i \cap \Sigma))_{< X}$.

In what follows, we will use the following version of Davenport's lemma
\cite{Davenport}.
\begin{prop}
\label{prop:Davenport}
Let $\cR$ be a bounded, semialgebraic multiset in $\R^n$ having maximum
multiplicity $m$ and that is defined by at most $k$ polynomial inequalities,
each having degree at most $l$. Let $L$ be a rank $n$ lattice inside $\R^n$. Then,
\[
\# (\cR \cap L) = \vol_L(\cR) + O(\max(\{\vol(\overline{\cR}),1\})),
\]
Here, $\vol_L$ is a constant multiple of $\vol$ with $\vol_L(\R^n/L) = 1$,
and $\vol(\overline{\cR})$ denotes the greatest $d$-dimensional volume of
any projection of $\cR$ onto a coordinate subspace obtained by equating
$n-d$ coordinates to zero, and where $d$ takes any value between $1$ and
$n-1$. The implied constant in the second summand depends only on $n,m,k$
and $l$.
\end{prop}

We now want to prove a ``cutting-off-the-cusp'' result in the style of
e.g. \cite[Proposition 8.11]{LagaThesis}, which should say that most
elements ``in the cusp'' fall into the subspace $W_0$ of $V$ defined in Section \ref{section:construct}. However,
unlike previous instances of this result, in our case we have multiple
cusps to worry about. Given that $\cS_1 \subset \cF \subset \cS_2$, we
can write $\cF = \cup_{j=1}^m \alpha_j \cS_j$, where the $\alpha_j \in G(F)$
are as in Section \ref{subs:FDomains} and $\cS_j$ are subsets of $\omega_jA_c'K_1$.
We can further assume that $\alpha_i \cS_i \cap \alpha_j \cS_j = \emptyset$
for all $i \neq j$. Then, we can write
\[
N(\Gamma,V_i,X) = \frac{1}{r_i \vol(G_0)} \sum_{j=1}^m \int_{g \in \alpha_j \cS_j}\# \{v \in V_0^{red} \cap (g\cB_X)\} dg
\]
For each cusp corresponding to $\alpha_j$, we can consider the weights of
the action of $\alpha_jT\alpha_j^{-1}$ on $V$ (the weight spaces
will be of the form $\alpha \cdot V_{\lambda}$, where $V_{\lambda}$
are the weight spaces for the action of $T$).
\begin{prop}
\label{prop:cusp}
Let $v_{0,j}$ denote the highest weight of $V$ under the action of $\alpha_jT\alpha_j^{-1}$. 
Then, there exists a constant $\delta > 0$ such that
\[
\int_{g \in \alpha_j \cS_j} \#\{v \in (V_0 \setminus (\alpha_j \cdot W_0)) \cap g \cB_X \mid v_{0,j}=0 \} dg = O(X^{\dim V - \delta})
\]
\end{prop}
\begin{proof}
This proof will follow the argument in \cite[\S 3.2]{BSWcoregular} and \cite[\S 3.3 and \S 5]{ThorneE6}.
Recall that $\cS_2$ is a finite cover of $\cF$ of absolutely bounded degree by the Siegel property (i.e. Proposition \ref{prop:Siegel}).
Hence, we can assume that
\[
\cS_j =  \omega_j A_c' K_1,
\]
following the notation in Section \ref{subs:FDomains}. 
Without loss of generality, we may assume that $\alpha_j = 1$,
since the statement for the rest of the cusps is analogous.
There exists a compact subset $\omega' \subset N(F_{\infty})$
that contains the union of all $t^{-1}\omega t$ as $t$ varies in $A_c'$.
Therefore, we have
\[
\int_{g \in \cS_j} \#\{v \in (V_0 \setminus W_0) \cap g \cB_X \mid v_0=0 \} dg \ll \int_{t \in A_c'} \#\{v \in (V_0 \setminus W_0) \cap t\omega'\cB_X \mid v_0=0 \} |\delta^{-1}(t)|dt
\]
Let $\Phi_V$ denote the characters of $V$ under the action of $T$. For two disjoint subsets
$M_0,M_1 \subset \Phi_V$, we define $S(M_0,M_1) = \{v \in V_0 \mid v_a = 0, \, \forall a \in M_0,\,v_b \neq 0, \, \forall b \in M_1\}$.
We define $\cC$ to be the collection of subsets $M_0 \subset \Phi_V$ such that if 
$a \in M_0$ and $b \geq a$, then $b \in M_0$. Additionally, given $M_0 \in \cC$,
we define $\lambda(M_0) = \{a \in \Phi_V \setminus M_0 \mid M_0 \cup \{a\} \in \cC\}$.
For $M_0 \in \cC$, we refer to a pair $(M_0,\lambda(M_0))$ as a cusp datum.

Any element in $(V_0 \setminus W_0) \cap g \cB_X$ falls inside one of
the subsets $S(M_0,M_1)$ for some cusp datum $(M_0,M_1)$. Therefore, it suffices to prove that
\[
N(\Gamma, S(M_0,M_1), X) = O(X^{\dim V - \delta})
\]
for every cusp datum such that $S(M_0,M_1) \nsubseteq W_0$.

Now, fix such a cusp datum $(M_0,M_1)$. We note that if $S(M_0,M_1) \cap g \cB_X$
is non-empty, then we must have $X|a(t)| \gg 1$ for every $a \in M_1$.
We also note that $\prod_{a \in \Phi_V} |a(t)| = 1$ for all $t \in T(F_{\infty})$.
Given that $V_0$ is commensurable with $V(\oO_K)$, there exists a constant
$C_0$ such that for every $\chi \in U_0$ and $v \in V_0$, we have that
either $v_{\chi} = 0$ or $|v_{\chi}| \geq C_0$. We define
\[
V_{M_0,M_1} := \{v \in V(F_{\infty}) \mid v_a = 0, \forall a \in M_0;\; |v_b| \geq C_0, \, \forall b \in M_1\}.
\]
Then, we have the following estimate:
\[
\vol(g \cB_X \cap V_{M_0,M_1}) = X^{\dim V - \# M_0} \prod_{a \in M_0}|a(t)|^{-1}.
\]
We also define
\[
A(M_0,M_1,X) := \{t \in T^{\theta}(F_{\infty}) \mid t \in A_c'; \; |a(t)| \gg X^{-1}, \forall a \in M_1\}.
\]
Finally, recall that $\delta^{-1}(t) = \prod_{a \in \Phi_G^+}a(t)$.
Using Proposition \ref{prop:Davenport} and the above observations, we get
\[
N(\Gamma, S(M_0,M_1), X) \ll X^{\dim V - \# M_0} \int_{t \in A(M_0,M_1,X)} \prod_{a \in \Phi_G^+}|a(t)|\prod_{a \in M_0}|a(t)|^{-1} dt.
\]
Therefore, we have reduced our statement to showing that
\[
\int_{t \in A(M_0,M_1,X)} \prod_{a \in \Phi_G^+}|a(t)|\prod_{a \in M_0}|a(t)|^{-1} dt = O(X^{\# M_0 - \delta}).
\]
Denote $w(t) = \prod_{a \in \Phi_G^+}|a(t)|\prod_{a \in M_0}|a(t)|^{-1}$.
Using a trick, due to Bhargava (cf. \cite[Lemma 19]{BhaQuintic}), 
let us consider a function $f \colon M_1 \to \R_{\geq 0}$.
We have that $\prod_{a \in M_1}(X|a(t)|)^{f(a)} \gg 1$, and therefore that
\begin{equation}
\label{eq:trick}
\int_{t \in A(M_0,M_1,X)} w(t) dt \ll X^{\sum_{a\in M_1} f(a)}\int_{t \in A(M_0,M_1,X)} w(t) \prod_{a \in M_1}|a(t)|^{f(a)} dt.
\end{equation}
Recall that any element $a \in X^{*}(T) \otimes \Q$ can be written
uniquely as $a = \sum_{\alpha_i \in S_G} n_i(a) \alpha_i$ for some rational
numbers $n_i(a)$. If we all the exponents of $w(t) \prod_{a \in M_1}|a(t)|^{f(a)}$
in terms of the basis $\alpha_i$ are positive, then the second integral
in \eqref{eq:trick} is bounded independently of $X$, and therefore
\[
\int_{t \in A(M_0,M_1,X)} w(t) dt \ll X^{\sum_{a\in M_1} f(a)}.
\]
Therefore, the proposition is reduced to finding a function $f \colon M_1 \to \R_{\geq 0}$
satisfying:
\begin{itemize}
\item We have $\sum_{a \in M_1} f(a) < \#M_0$.
\item For all $i$, we have $\sum_{\beta \in \Phi_G^+} n_i(\beta) - \sum_{a \in M_0} n_i(a) + \sum_{a \in M_1} f(a) n_i(a) > 0$.
\end{itemize}
Note that this last condition is independent of the base field $F$, so
it is sufficient to prove the cutting-off-the-cusp results over $\Q$.
This is the content of \cite[Proposition 8.21]{LagaThesis}.
\end{proof}
Hence, if an element is in the cusp (i.e. the highest weight has coefficient zero),
then it always falls in $W_0$ (and is therefore reducible) expect for
negligibly many cases. We will now see that, analogously, almost all
the elements in the main body are irreducible:
\begin{prop}
\label{prop:Selberg1}
There exists a constant $\delta > 0$ such that
\[
\int_{g \in \cF} \#\{v \in V_0^{red} \cap g \cB_X \mid v_0\neq0 \} dg = O(X^{\dim V - \delta}).
\]
\end{prop}
\begin{proof}
To prove this, we will use the Selberg sieve, in an analogous way to
\cite{STSelberg}.
The general argument in \cite{STSelberg}, for $F = \Q$, is the following: 
assume that for any translate $L$ of $mV(\Z)$ we have that
\begin{equation}
\label{eq:Selberg}
N^*(L \cap V_i, X) = c_i m^{-A}X^B + O(m^{-C}X^{B-D}),
\end{equation}
where $A,B,C,D$ and $c_i$ are positive constants, and $N^*$ is some orbit-counting
function. Let $S = \cap_{p \text{ prime}} S_p \subset V(\Z)$
be a set defined by infinitely many congruence conditions modulo $p$, with each set $S_p$ having
density $\lambda_p$. Assume that $\lambda_p$ tends to some constant $\lambda \in (0,1)$
as $p$ tends to infinity. Then, it is shown that
\[
N^*(S \cap V_i,X) = O(X^{B-\delta}),
\]
for some constant $\delta > 0$ which can be obtained explicitly depending
on $A,B,C,D$. If instead of working over $\Q$ we work over a number field
$F$, we can apply the same argument substituting our sets $S_p$ for
$S_{\pp}$ for $\pp$ a prime ideal of $\oO_{F}$, and using the version
of the Selberg sieve stated in \cite[Satz 1]{SelbergNF}.

In our case, we set
\[
N^*(S,X) = \int_{g \in \cF} \# \{v \in S \cap g\cB_X \mid v_0 \neq 0 \} dg.
\]
First, we need a power saving estimate for $N^*(L \cap V_i, X)$, where
$L$ is a translate of $I V_0$ for some ideal $I \subset \oO_F$. This is
done in Corollary \ref{cor:irred}. We now let
\[
S_{\pp} = \{v \in V(k_{\pp}) \mid \Delta(v) = 0 \text{ or } v \text{ is }k_{\pp}\text{-reducible}\}.
\]
By \cite[Theorem 7.16]{LagaThesis}, any reducible element in $V_0$ falls
into the set $S_{\pp}$ for all primes not dividing $N_{bad}$, and 
\cite[Proof of Lemma 8.20]{LagaThesis} shows that the density of the sets
$S_{\pp}$ tends to some constant $c \in (0,1)$ as $N\pp$ tends to infinity.
In conclusion, we can apply the Selberg sieve to obtain the desired power saving estimate.
\end{proof}

Finally, we argue that instead of integrating over the fundamental domain
$\cF$, it is enough to integrate over the smaller and more convenient set
$\cS_1$:
\begin{lemma}
There is some constant $\delta > 0$ such that
\[
\int_{g \in \cF\setminus \cS_1} \# \{v \in V_0^{red} \cap g \cB_{X} \} dg = O(X^{\dim V - \delta}).
\]
\end{lemma}
\begin{proof}
It suffices to do so in each cusp separately, so fix a cusp corresponding
to $\alpha_j$. The region of integration in this case is a subset of
\[
\alpha_j\omega_j\{t \in A_c' \mid |\alpha(t)|\geq c' \text{for some }\alpha \in S_G\}K_1,
\]
for some choice of $c,c'$. The computations in the proof of Proposition
\ref{prop:cusp} directly show that the integral in this case is $O(X^{\dim V-\delta})$
for some $\delta > 0$, as wanted.
\end{proof}

\subsection{Slicing}
\label{subs:slicing}

The results in Section \ref{subs:reductions} show that
\begin{equation}
\label{eq:slicing}
N(\Gamma,V_i,X) = \frac{1}{r_i \vol(G_0)} \sum_{j=1}^m\int_{g \in \alpha_j \omega_j A_c'K_1} \# \{v \in (\alpha_j\cdot W_0) \cap g \cB_{X} \} dg + O(X^{\dim V - \delta})
\end{equation}
for some constant $\delta > 0$. To estimate the number of lattice points
in the given region, we would want to use Proposition \ref{prop:Davenport}.
However, we can't use it directly, because some of the projections onto coordinate
hyperplanes of $W_0$ are of
the order of the main term! To circumvent this, we will slice the region
$W_0$ according to the values of the height-one coefficients.

We will compute the integral in \eqref{eq:slicing} separately for each
cusp $\alpha_j \omega_j A_c'K_1$. To simplify notation, we will work
with the cusp with $\alpha_j = 1$, since the computations for the rest
are analogous.

Let $v \in W_0(F_{\infty})$, and let $v_1,\dots,v_k$ denote its height-one
coefficients, where $k = \# S_G$. For any $b \in (F_{\infty})^k$ and any
subset $\cS \subset W_0(F_{\infty})$, we write
\[
\cS_b = \{v \in \cS \mid (v_1,\dots,v_k) = b\}.
\]
For any $v \in V_0 \cap W_0$, we know that the values of its height-one
coefficients fall into some lattice $\cL \subset F^k$ which as an additive
subgroup is commensurable to $\oO_F^k$. We can then write
\[
\#(g \cB_X \cap W_0 \cap V_0) = \sum_{b \in \cL}\#(g \cB_X \cap W_0 \cap V_0)_b
\]
In fact, we can avoid summing over some of the $b$:
\begin{lemma}
Let $v \in W_0(F_{\infty})$. If $v_i = 0$ for some height-one coefficient $i$, then $\Delta(v) = 0$.
\end{lemma}
\begin{proof}
It suffices to consider each completion $F_w$ separately, for all infinite
places $w$.
Let $\{\alpha_1,\dots,\alpha_k\}$ be the height-one weights, and assume
that the coefficient of $\alpha_i$ of $v$ is zero in $F_v$. Let 
$\lambda_i\colon \G_m \to G_{\C}$ be the one-parameter subgroup 
such that $(\alpha_j \circ \lambda_i)(t) = t^{\delta_{ij}}$. Then,
$v$ has no positive weights with respect to $\lambda_i$, and so by Proposition
\ref{prop:H-M} we get the result.
\end{proof}
We can now use Davenport's lemma, using the natural bijection between 
$F_{\infty}$ and $\R^{[F:\Q]}$, to estimate
\[
\#(g \cB_X \cap W_0 \cap V_0)_b = \vol((g\cB_X)|_b)(1+O(X^{-1/[F:\Q]})).
\]
The weight of any non-height-one coefficient in $W_0$ is $\gg 1$, so 
the range of values of any real coefficient varying in $\cB_X$ is 
$\gg X^{1/[F:\Q]}$, so we get an error term of the order of $O(X^{-1/[F:\Q]})$.

Now, note that $KG_0 = G_0$ and that unipotent transformations preserve
both the value of the height-one coefficients and the volume, so we get
that $\vol((ntk\cB_X)|_b) = \vol((t\cB_X)|_b)$, where $t \in A_c'$. It
will be convenient for us to integrate not over $A_c'$, but rather over
$T \cap \cS_1$, which is a set of the form $A_c' \times K_T$, for some
subset $K_T \subset K_1$. We note that by construction of $K_1$, the set
$A_c'\times K_T$ corresponds to a fundamental domain for the action of
$\Gamma \cap T$ on $T(F_{\infty})$ (cf. Section \ref{subs:corners}). There is
a natural measure on $A_c'\times K_T$, inherited by restriction of the
measures $dt$ and $dk$: we will denote it $dt$ by abuse of notation.
Then, given that $\omega_j$ and $K_1$ have finite measure, we get that
\begin{equation}
\label{eq:slicing1}
N(\Gamma,V_i,X) = C\sum_{\substack{b \in \cL \\ b_i \neq 0\, \forall i}} \int_{t \in A_c'\times K_T} \vol((t\cB_X)|_b) |\delta(t)|^{-1} dt,
\end{equation}
for some constant $C$.
For each height-one coefficient $v_i$, its weight under the
action of $T$ is $\alpha_i(t)$. We will let $\beta_i := b_i/(X^{1/[F:\Q]} \alpha_i(t))$,
and $\beta = (\beta_i)_i$. Denote by $\Phi_V$ the different weights of
the action of $T$ on $V$, and by $\Phi_V^-$ the negative weights. Then,
we have that
\[
\vol((t\cB_X)|_b) = \vol(t X^{1/[F:\Q]}\cdot(\cB_1)|_{\beta}) = X^{\dim W_{\flat}} \prod_{\alpha \in \Phi_V^-}|\alpha(t)| \vol(\cB|_{\beta}).
\] 
We will make the change of variables $t \mapsto \beta = (\beta_1,\dots,\beta_k)$, 
under which $dt = d \beta$, where $d \beta = \prod_{i=1}^k \frac{d\beta_i}{|\beta_i|}$. 
In Section \ref{subs:cases}, we will 
explicitly compute the volume of the cuspidal region for each of the 
possible cases. We will obtain a polynomial $Z(\beta) = \prod_{i}\beta_i^{e_i}$
with integer exponents $e_i \geq 2$, and we will see that
\begin{equation}
\label{eq:explicit}
X^{\#\Phi_V^-} \prod_{\alpha \in \Phi_V^-}|\alpha(t)| |\delta^{-1}(t)| = X^{\dim V} \frac{|Z(\beta)|}{|Z(b)|}.
\end{equation}
The group $\Gamma_T := \Gamma \cap T$ acts naturally on $\cL$ by 
$t\cdot (b_1,\dots,b_k) = (\alpha_1(t)b_1,\dots,\alpha_k(t)b_k)$ for
$t \in \Gamma_T$ and $b \in \cL$ (we know that $t \cdot b \in \cL$ because
$V_0$ is invariant under $\Gamma$).
Denote $A' = \cup_{c > 0}A_c'$ (i.e. the region defined by $A_c'$ without
the condition that $|\alpha(t)|\leq c$).
The change of variables from $t$ to $\beta$ sends the domain of integration
$A' \times K_T$ to some region $Y_{b} \subset F_{\infty}^{|M_{\infty}|}$,
and say that the region $A_{c}'\times K_T$ gets sent to $Y_{b} \setminus Y_c$.
It follows that
\[
\int_{t \in  A_c'\times K_T} \vol((t\cB_X)|_b)  |\delta^{-1}(t)| dt = \frac{X^{\dim V}}{|Z(b)|} \int_{\beta \in Y_{b} \setminus Y_c} |Z(\beta)| \vol(\cB|_{\beta}) d\beta.
\]
The set $Y_c$ corresponds to elements $\beta$ with $|\beta| \gg X^{-1}$
-- in particular, the integral over $Y_c$ is bounded by $O(X^{-1})$, so
it can be added to the error term. Now, for any $b_0 \in \cL$, we have
that at most $O_{\eps}(X^{\eps})$ choices of $b \in \Gamma_t b_0$ give
a non-zero volume of $\cB|_{\beta}$. Therefore, we can write
\[
\sum_{b \in \Gamma_T b_0} \int_{t \in  A_c'\times K_T} \vol((t\cB_X)|_b)  |\delta^{-1}(t)| dt = \sum_{b \in \Gamma_T b_0}\frac{X^{\dim V}}{|Z(b)|} \int_{\beta \in Y_{b}} |Z(\beta)| \vol(\cB|_{\beta}) d\beta + O_{\eps}(X^{-1+\eps}).
\]
Now, $\cup_{b \in \Gamma_T b_0} Y_{b} = (F_{\infty}^{|M_{\infty}|})^k$, so we
get:
\begin{align*}
\sum_{b \in \Gamma_T b_0} \int_{t \in  A_c'\times K_T} \vol((t\cB_X)|_b)  |\delta^{-1}(t)| dt &=\frac{X^{\dim V}}{|Z(b_0)|} \int_{\beta \in (F_{\infty}^{|M_{\infty}|})^k} |Z(\beta)| \vol(\cB|_{\beta}) d\beta + O_{\eps}(X^{-1+\eps}) \\
&= \frac{X^{\dim V}}{|Z(b_0)|} \int_{v \in W_{0}} |Z^{\times}(v)| dv + O_{\eps}(X^{-1+\eps}). \\
\end{align*}
Here, $Z^{\times}(\beta) = \prod_{i=1}^k \beta_i^{e_i-1}$, where the
$e_i$ are the exponents appearing in $Z(\beta)$. Adding
over all $b_0 \in \cL/\Gamma_T$ and combining with \eqref{eq:slicing1},
we conclude the proof of Theorem \ref{theo:constant}.

\subsection{Congruence conditions}

In the sequel, it will be convenient for us to not only have an estimate
for the number of reducible orbits, but we will also need some knowledge
about what happens when we impose finitely many congruence conditions
in our orbits. It will suffice to do our analysis in the cusp; to that
effect, consider the counting function
\[
N^{cusp}(\Gamma,V_0,X) = \int_{g \in \cF} \#\{v \in V_0 \cap W_0 \cap (g\cB_X)\} dg.
\]
In the previous sections, we obtained that
\[
N^{cusp}(\Gamma,V_0,X) = CX^{\dim V} + O(X^{\dim V - \frac{1}{[F:\Q]}}),
\]
where $C$ and the implied constant depend only on the choice of $\Gamma$
and $V_0$. We will now obtain the following:
\begin{theorem}
Let $L$ be a translate of $IV_0$, for some ideal $I$ of $\oO_F$. If
\[
N^{cusp}(\Gamma,V_0,X) = CX^{\dim V} + O(X^{\dim V - \frac{1}{[F:\Q]}}),
\]
then
\[
N^{cusp}(\Gamma,L,X) = C(NI)^{-\dim V}X^{\dim V} + O((NI)^{\frac{1}{[F:\Q]}-\dim V}X^{\dim V - \frac{1}{[F:\Q]}}).
\]
The implied constant is independent of the choice of $L$.
\end{theorem}
\begin{proof}
The argument goes through in the same way as Section \ref{subs:slicing},
and the only difference is in the application of Davenport's lemma, where
the additional terms appear by taking care of the change of lattices.
\end{proof}
 
\subsection{Case-by-case analysis}
\label{subs:cases}

In this section, we complete the proof of Theorem \ref{theo:constant} by performing a case-by-case analysis.
For the $D_n$ and $E_n$ cases, we will explicitly compute the dimension
and volume of $W_{\flat}$ (which was defined to be the set of coefficients
of $W_0$ of non-positive height), and the modular function
$\delta(t) = \prod_{\beta \in \Phi_G^{-}} \beta(t) = \det \Ad(t)|_{\Lie N(F_{\infty})}$.

\subsubsection{$D_{2n+1}$}

The exposition in the $D_n$ cases is inspired by \cite[Appendix A]{LagaThesis}
and \cite[\S 7.2.1]{ShankarD2n+1}.
We start by describing explicitly the representation $(G,V)$ of $D_{2n+1}$
in the form given by Table \ref{table:explicit}.

Let $n \geq 2$ be an integer. Let $U_1$ be a $\Q$-vector space with basis
$\{e_1,\dots,e_n,u_1,e_n^{*},\dots,e_1^{*}\}$, endowed with the symmetric
bilinear form $b_1$ satisfying $b_1(e_i,e_j) = b_1(e_i,u_1) = b_1(e_i^*,e_j^*) = b_1(e_i^*,u_1) = 0$,
$b_1(e_i,e_j^*) = \delta_{ij}$ and $b_1(u_1,u_1) = 1$ for all $1 \leq i,j \leq n$.
In this case, given a linear map $A \colon U \to U$ we can define its \emph{adjoint}
as the unique map $A^{\ast} \colon U \to U$ satisfying $b_1(Av,w) = b_1(v,A^{\ast}w)$
for all $v,w \in U$. In terms of matrices, $A^*$ corresponds to taking 
the reflection of $A$ along its antidiagonal when working with the fixed basis. We can define 
$\So(U_1,b_1) := \{ g \in \SL(U_1) \mid gg^* = \id\}$, with a Lie algebra
that can be identified with $\{A \in \End(U) \mid A = -A^* \}$.

Let $U_2$ be a $\Q$-vector space with basis $\{f_1,\dots,f_n,u_2,f_n^{\ast},\dots,f_1^{\ast}\}$,
with a similarly defined bilinear form $b_2$. Let $(U,b) = (U_1,b_1) \oplus (U_2,b_2)$.
Let $H = \So(U,b)$, and consider $\hh := \Lie H$. With respect to the basis
\[
\{e_1,\dots,e_n,u_1,e_n^{\ast},\dots,e_1^{*},f_1,\dots,f_n,u_2,f_n^{\ast},\dots,f_1^{\ast}\},
\]
the adjoint of a block matrix according to the bilinear form $b$ is given by
\[
\begin{pmatrix}
A & B \\ C & D
\end{pmatrix}^*
=
\begin{pmatrix}
A^* & C^* \\ B^* & D^*
\end{pmatrix},
\]
where $A^*,B^*,C^*,D^*$ denote reflection by the antidiagonal. An element
of $\hh$ is given by
\[
\left\{\begin{pmatrix}
B & A \\ -A^* & C
\end{pmatrix} \,\middle\vert\; B = -B^*, \, C = -C^*\right\}.
\]
The stable involution $\theta$ is given
by conjugation by $\diag(1,\dots,1,-1,\dots,-1)$, where the first $2n+1$
entries are $1$ and the last $2n+1$ entries are given by $-1$. Under this
description, we see that
\[
V = \left\{\begin{pmatrix}
0 & A \\ -A^* & 0
\end{pmatrix} \,\middle\vert\; A \in \text{Mat}_{(2n+1)\times(2n+1)}\right\}.
\]
Moreover, $G = (H^{\theta})^{\circ}$ is isomorphic to $\So(U_1) \times \So(U_2)$.
We will use the map
\[
\begin{pmatrix}
0 & A \\ -A^* & 0
\end{pmatrix} \mapsto A
\]
to establish a bijection between $V$ and $\Hom(U_2,U_1)$, where $(g,h) \in \So(U_1) \times \So(U_2)$
acts on $A \in V$ as $(g,h)\cdot A = gAh^{-1}$.

Let $T$ be the maximal torus $\diag(t_1,\dots,t_n,1,t_n^{-1},\dots,t_1^{-1},s_1,\dots,s_n,1,s_n^{-1},\dots,s_1^{-1})$
of $G$. A basis of simple roots for $G$ is
\[
S_G = \{t_1-t_2,\dots,t_{n-1}-t_n\} \cup \{s_1-s_2,\dots,s_{n-1}-s_n\}.
\]
A positive root basis for $V$ can be taken to be
\[
S_V = \{t_1-s_1,s_1-t_2,\dots,t_n-s_n,s_n\}.
\]
For convenience, we now switch to multiplicative notation for the roots.
We make the change of variables $\alpha_i = t_i/t_{i+1}$ for $i = 1,\dots,n-1$
and $\alpha_n = t_n$; similarly $\gamma_i = s_i/s_{i+1}$ for $i = 1,\dots,n-1$
and $\gamma_n = s_n$. The estimate for the volume becomes
\[
\prod_{\lambda \in \Phi_V^-} X|\lambda(t)| = X^{2n^2+2n+1} \prod_{i=1}^n |\alpha_i|^{-2in+i^2-2i}|\gamma_i|^{-2in+i^2}.
\]
The modular function in our case is
\[
|\delta^{-1}(t)| = \prod_{i=1}^n |\alpha_i|^{2in-i^2}|\gamma_i|^{2in-i^2}.
\]
Changing variables to $\beta_i = b_i/(X \lambda_i(v_i))$, where $\lambda_i$ are the 
height-one weights, we obtain
\[
\prod_{\lambda \in \Phi_V^-} X|\lambda(t)| \delta^{-1}(t) = X^{4n^2+4n+1} \frac{|Z(\beta)|}{|Z(b)|},
\]
where $Z(\beta) := \prod_{i=1}^n (\beta_{2i-1}\beta_{2i})^{2i}$.

\subsubsection{$D_{2n}$}
The analysis in this case is very similar to the $D_{2n+1}$ case. Now, we
consider the $\Q$-vector space $U_1$ with basis $\{e_1,\dots,e_n,e_n^*,\dots,e_1^*\}$,
endowed with a symmetric bilinear form $b_1(e_i,e_j) = b_1(e_i^*,e_j^*) = 0$, 
$b_1(e_i,e_j^*) = \delta_{ij}$. We also consider a $\Q$-vector space $U_2$
with basis $\{f_1,\dots,f_n,f_n^*,\dots,f_1^*\}$, with an analogous symmetric
bilinear form $b_2$.

Let $(U,b) = (U_1,b_1) \oplus (U_2,b_2)$, let $H' = \So(U,b)$ and define
$H$ to be the quotient of $H'$ by its centre of order 2. Under the basis
\[
\{e_1,\dots,e_n,e_n^*,\dots,e_1^*,f_1,\dots,f_n,f_n^*,\dots,f_1^*\},
\]
the stable involution is given by conjugation with $\diag(1,\dots,1,-1,\dots,-1)$.
Similarly to the $D_{2n+1}$ case, we have
\[
V = \left\{\begin{pmatrix}
0 & A \\ -A^* & 0
\end{pmatrix} \,\middle|\; A \in \text{Mat}_{2n \times 2n}\right\},
\]
where $A^*$ denotes reflection by the antidiagonal. In this case, the
group $G = (H^{\theta})^{\circ}$ is isomorphic to $\So(U_1) \times \So(U_2)/\Delta(\mu_2)$,
where $\Delta(\mu_2)$ denotes the diagonal inclusion of $\mu_2$ into the
centre $\mu_2 \times \mu_2$ of $\So(U_1) \times \So(U_2)$. As before,
we can identify $V$ with the space of $2n\times2n$ matrices using the
map
\[
\begin{pmatrix}
0 & A \\ -A^* & 0
\end{pmatrix} \mapsto A,
\]
where $(g,h) \in G$ acts by $(g,h)\cdot A = gAh^{-1}$.

We consider the maximal torus $T$ of $H$ given by $\diag(t_1,\dots,t_n,t_n^{-1},\dots,t_1^{-1},s_1,\dots,s_n,s_n^{-1},\dots,s_1^{-1})$.
A basis of simple roots for $H$ and $G$ are given by
\begin{align*}
&S_H = \{t_1-s_1,s_1-t_2,\dots,s_{n-1}-t_n,t_n-s_n,s_n+t_n\}, \\
&S_G = \{t_1-t_2,\dots,t_{n-1}-t_n,t_{n-1}+t_n\} \cup \{s_1-s_2,\dots,s_{n-1}-s_n,s_{n-1}+s_n\}.
\end{align*}

Let $\alpha_i = t_i/t_{i+1}$ and $\gamma_i = s_i/s_{i+1}$ for $i = 1,\dots,n$,
and let $\alpha_n = t_{n-1}t_n$ and $\gamma_n = s_{n-1}s_n$. Under this change
of variables, the estimate for the volume is:
\[
\prod_{\lambda \in \Phi_V^-}X|\lambda(t)| = X^{2n^2} \left(\prod_{i=1}^{n-2}|\alpha_i|^{-2in+i^2-i}|\alpha_{n-1}|^{(-n^2-n+4)/2}|\alpha_n|^{(-n^2-n)/2}\prod_{i=1}^{n-2}|\gamma_i|^{-2in+i^2+i}|\gamma_{n-1}\gamma_n|^{(-n^2+n)/2}\right).
\]
The modular function is
\[
|\delta^{-1}(t)| = \prod_{i=1}^{n-2} |\alpha_i\gamma_i|^{i^2-2in+i} |\alpha_{n-1}\gamma_{n-1}\alpha_n\gamma_n|^{-(n-1)n/2}.
\]
As before, we can compute:
\[
\prod_{\lambda \in \Phi_V^-}X|\lambda(t)| |\delta^{-1}(t)| = X^{4n^2} \frac{|Z(\beta)|}{|Z(b)|},
\]
where $Z(\beta) = \prod_{i=1}^{n-1}(\beta_{2i-1}\beta_{2i})^{2i}\cdot(\beta_{2n-1}\beta_{2n})^n$.
\subsubsection{$E_6$}
For the $E_6$ case, we use the conventions and computations in \cite[\S 2.3, \S 5]{ThorneE6}.

Let $S_H = \{\alpha_1,\dots,\alpha_6\}$, where the Dynkin diagram of $H$ is:
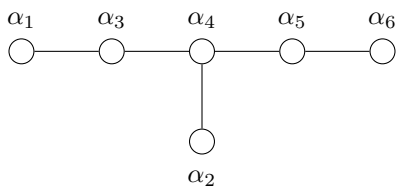
\begin{center}
\begin{tikzpicture}[scale=0.6]
  \node[circle, draw,label={$\alpha_1$}] (A1) at (0,0) {};
  \node[circle, draw,label={$\alpha_3$}] (A3) at (2,0) {};
  \node[circle, draw,label={$\alpha_4$}] (A4) at (4,0) {};
  \node[circle, draw,label={$\alpha_5$}] (A5) at (6,0) {};
  \node[circle, draw,label={$\alpha_6$}] (A6) at (8,0) {};
  \node[circle, draw,label=below:{$\alpha_2$}] (A2) at (4,-2) {};
  
  \draw[-] (A1) -- (A3) -- (A4) -- (A5) -- (A6);
  \draw[-] (A2) -- (A4);
\end{tikzpicture}
\end{center}
The pinned automorphism $\vartheta$ consists of a reflection
around the vertical axis. We can define a root basis $S_G = \{\gamma_1,\gamma_2,\gamma_3,\gamma_4\}$
of $G$ as $\gamma_1 = \alpha_3+\alpha_4$, $\gamma_2 = \alpha_1$, $\gamma_3 = \alpha_3$ and $\gamma_4 = \alpha_2 + \alpha_4$.
Under this basis, we have
\[
\prod_{\lambda \in \Phi_V^-}X|\lambda(t)| = X^{22} |\gamma_1|^{-12}|\gamma_2|^{-18}|\gamma_3|^{-22}|\gamma_4|^{-12}
\]
The modular function is
\[
|\delta^{-1}(t)| = |\gamma_1|^{8}|\gamma_2|^{14}|\gamma_3|^{18}|\gamma_4|^{10}.
\]
The weights of the height-one coefficients are $\{\gamma_2,-\gamma_1+\gamma_3+\gamma_4,\gamma_3,\gamma_1-\gamma_3\}$.
In light of this, we obtain
\[
\prod_{\lambda \in \Phi_V^-}X|\lambda(t)||\delta^{-1}(t)| = X^{42}\frac{|Z(\beta)|}{|Z(b)|}.
\]
where $Z(\beta) = \beta_1^4 \beta_2^2 \beta_3^8 \beta_4^6$.
\subsubsection{$E_7$}
For the $E_7$ and $E_8$ cases, we follow the conventions in \cite{RTE78}.
Let $S_H = \{\alpha_1,\dots,\alpha_7\}$, where the Dynkin diagram of $H$ is:
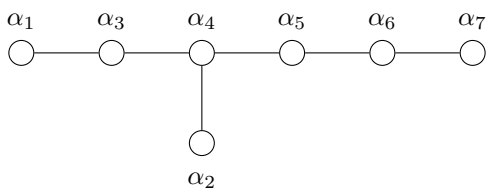
\begin{center}
\begin{tikzpicture}[scale=0.6]
  \node[circle, draw,label={$\alpha_1$}] (A1) at (0,0) {};
  \node[circle, draw,label={$\alpha_3$}] (A3) at (2,0) {};
  \node[circle, draw,label={$\alpha_4$}] (A4) at (4,0) {};
  \node[circle, draw,label={$\alpha_5$}] (A5) at (6,0) {};
  \node[circle, draw,label={$\alpha_6$}] (A6) at (8,0) {};
  \node[circle, draw,label={$\alpha_7$}] (A7) at (10,0) {};
  \node[circle, draw,label=below:{$\alpha_2$}] (A2) at (4,-2) {};
  
  \draw[-] (A1) -- (A3) -- (A4) -- (A5) -- (A6) -- (A7);
  \draw[-] (A2) -- (A4);
\end{tikzpicture}
\end{center}
The root basis $S_G = \{\gamma_1,\dots,\gamma_7\}$ can be described as
\begin{align*}
\gamma_1 &= \alpha_3 + \alpha_4 \\
\gamma_2 &= \alpha_5 + \alpha_6 \\
\gamma_3 &= \alpha_2 + \alpha_4 \\
\gamma_4 &= \alpha_1 + \alpha_3 \\
\gamma_5 &= \alpha_4 + \alpha_5 \\
\gamma_6 &= \alpha_6 + \alpha_7 \\
\gamma_7 &= \alpha_2 + \alpha_3 + \alpha_4 + \alpha_5 \\
\end{align*}
The volume of $W_{\flat}$ can be computed to be
\[
\prod_{\lambda \in \Phi_V^-}X|\lambda(t)|= X^{35} |\gamma_1|^{-15/2}|\gamma_2|^{-13}|\gamma_3|^{-33/2}|\gamma_4|^{-18}|\gamma_5|^{-35/2}|\gamma_6|^{-15}|\gamma_7|^{-21/2}.
\]
The modular function for $G$ can be computed to be
\[
\delta^{-1}(t) = |\gamma_1|^7|\gamma_2|^{12}|\gamma_3|^{15}|\gamma_4|^{16}|\gamma_5|^{15}|\gamma_6|^{12}|\gamma_7|^7.
\]
We can compute the weights $\beta_i$ corresponding to the height-one
coefficients, with the end result being
\[
\prod_{\lambda \in \Phi_V^-}X|\lambda(t)||\delta^{-1}(t)| = X^{70}\frac{|Z(\beta)|}{|Z(b)|},
\] 
for $Z(\beta) = \beta_1^2 \beta_2^5 \beta_3^6 \beta_4^8 \beta_5^7 \beta_6^4 \beta_7^3$.
\subsubsection{$E_8$}
Let $S_H = \{\alpha_1,\dots,\alpha_8\}$, where the Dynkin diagram of $H$ is:
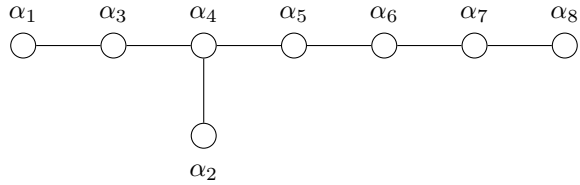
\begin{center}
\begin{tikzpicture}[scale=0.6]
  \node[circle, draw,label={$\alpha_1$}] (A1) at (0,0) {};
  \node[circle, draw,label={$\alpha_3$}] (A3) at (2,0) {};
  \node[circle, draw,label={$\alpha_4$}] (A4) at (4,0) {};
  \node[circle, draw,label={$\alpha_5$}] (A5) at (6,0) {};
  \node[circle, draw,label={$\alpha_6$}] (A6) at (8,0) {};
  \node[circle, draw,label={$\alpha_7$}] (A7) at (10,0) {};
  \node[circle, draw,label={$\alpha_8$}] (A8) at (12,0) {};
  \node[circle, draw,label=below:{$\alpha_2$}] (A2) at (4,-2) {};
  
  \draw[-] (A1) -- (A3) -- (A4) -- (A5) -- (A6) -- (A7) -- (A8);
  \draw[-] (A2) -- (A4);
\end{tikzpicture}
\end{center}
The root basis $S_G = \{\gamma_1,\dots,\gamma_8\}$ can be described as
\begin{align*}
\gamma_1 &= \alpha_2 + \alpha_3 + \alpha_4 + \alpha_5 \\
\gamma_2 &= \alpha_6 + \alpha_7 \\
\gamma_3 &= \alpha_4 + \alpha_5 \\
\gamma_4 &= \alpha_1 + \alpha_3 \\
\gamma_5 &= \alpha_2 + \alpha_4 \\
\gamma_6 &= \alpha_5 + \alpha_6 \\
\gamma_7 &= \alpha_7 + \alpha_8 \\
\gamma_8 &= \alpha_3 + \alpha_4\\
\end{align*}
The volume of $W_{\flat}$ can be computed to be
\[
\prod_{\lambda \in \Phi_V^-}X|\lambda(t)| = X^{64} |\gamma_1|^{-18}|\gamma_2|^{-30}|\gamma_3|^{-40}|\gamma_4|^{-48}|\gamma_5|^{-54}|\gamma_6|^{-58}|\gamma_7|^{-30}|\gamma_8|^{-30}.
\]
The modular function for $G$ can be computed to be
\[
|\delta^{-1}(t)| = |\gamma_1|^{14}|\gamma_2|^{26}|\gamma_3|^{36}|\gamma_4|^{44}|\gamma_5|^{50}|\gamma_6|^{54}|\gamma_7|^{28}|\gamma_8|^{28}.
\]
We get
\[
\prod_{\lambda \in \Phi_V^-}X|\lambda(t)||\delta^{-1}(t)| = X^{128}\frac{|Z(\beta)|}{|Z(b)|},
\]
with $Z(\beta) = \beta_1^4 \beta_2^8 \beta_3^{10} \beta_4^{14} \beta_5^{12} \beta_6^8 \beta_7^6 \beta_8^2$.

\section{Proof of the main results}
\label{section:final}

We are now in a position to prove the main results.

\subsection{Elements with big stabiliser}
We first proof a necessary result about elements with a big stabiliser.
As in Section \ref{section:Counting}, let $V_0$ be a lattice inside $V(F_{\infty})$ which
is commensurable with $V(\oO_K)$, and let $\Gamma$ be an arithmetic subgroup
of $G(F)$ which preserves $V_0$. We denote by $V_0^{bs}$ the set of elements
$v \in V_0$ which satisfy $\#\Stab_{G(F)}v > 1$.
\begin{prop}
\label{prop:bs}
There exists a constant $\delta_{bs} > 0$ with
\[
N(\Gamma,V_0^{bs,red},X) = O(X^{\dim V - \delta_{bs}}).
\]
\end{prop}
\begin{proof}
We can see that the density of elements in $V(\oO_{\pp})$ having big 
stabiliser tends to some constant in $(0,1)$ by the same argument as in 
\cite[Proof of Lemma 8.20]{LagaThesis}. Then, we can use the Selberg sieve
as in Proposition \ref{prop:Selberg1}, now using the estimate in
Theorem \ref{theo:constant}.
\end{proof}
\begin{obs}
We remark that we could not have proven Proposition \ref{prop:bs} at
the same time as Proposition \ref{prop:Selberg1}, as we need the estimate
in Theorem \ref{theo:constant} to apply the Selberg sieve in this case.
\end{obs}

\subsection{Tail estimates}
\label{subs:tail}
To prove Theorem \ref{theorem:tail}, we need to obtain tail estimates for
the sets $\cW_I^{(1)}$ and $\cW_I^{(2)}$. The required estimate for
$\cW_I^{(1)}$ can be obtained using \cite[Theorem 18]{BSWprehom} (which is
essentially by following the argument in \cite[Theorem 3.3]{BEkedahl}).
For $\cW_I^{(2)}$, recall that we only deal with ideals $I$ which are
coprime to a certain element $N_{bad}$, as explained in Section \ref{section:construct}.
We can prove the following:
\begin{theorem}
There exists a constant $\delta > 0$ such that
\[
\sum_{\substack{I \text{ squarefree}\\NI > M \\ (I,N_{bad}) = 1}} \#\{b \in \G_m(F) \backslash\cW_I^{(2)} \mid \Ht(b) < X\} = O\lp\frac{X^{\dim V}}{M}\rp + O(X^{\dim V - \delta}).
\]
\end{theorem}
\begin{proof}
In Section \ref{section:construct}, for every $\beta_i \in \cl(G)$ we constructed
\[
W_{i,M} := \left\{v \in V_{\beta_i}\,\middle|\; v = g_I\kappa_b,\, I \text{ squarefree},\, I \text{ coprime to }N_{bad}, \, NI > M,\, g_I \in G_I,\, b \in B(\oO_K)\right\}.
\]
By the results in Section \ref{section:construct}, it is sufficient to
obtain an appropriate bound on the number of $G_{\beta_i}$-orbits of $W_{i,M}$ with 
invariants in $\Sigma$.
By following the same averaging argument as in Section \ref{section:Counting},
we get that
\[
N(G_{\beta_i},W_{i,M},X) = \int_{g \in \cF} \#\{v \in W_{i,M} \cap g\cB_X\} dg.
\]
We can carry out the same argument, now assuming that we can restrict to
those elements in $W_{i,M}$ with trivial stabiliser by Proposition \ref{prop:bs}.
If $v \in W_{i,M} \cap W_0$ has trivial stabiliser, then $Q(v) \gg M$ by
Proposition \ref{prop:p-orbit}, or alternatively $Z(v) \gg M^2$. We have
that:
\[
N(G_{\beta_i},W_{i,M},X) \ll \sum_{\substack{b_0 \in \cL/\Gamma_T\\Z(b_0) \gg M^2}} \frac{1}{|Z(b_0)|} + O(X^{\dim V - \delta}),
\]
and the written sum is $O(1/M)$, which concludes the proof.
\end{proof}
Thus, we have concluded the proof of Theorem \ref{theorem:tail}.
We can combine both estimates
for the strongly divisible case and the weakly divisible case. For $I$
a squarefree ideal, denote by $\cW_I$ the set of elements $b \in B(\oO_F)$
such that $I^2$ divides $\Delta(b)$. Then, in the style of \cite[Theorem 4.4]{BSWsquarefree},
we obtain:
\begin{theorem}
\label{theo:fullTail}
There is a constant $\delta > 0$ such that
\[
\sum_{\substack{I \text{ squarefree}\\NI > M \\ (I,N_{bad}) = 1}} \#\{b \in \G_m(F) \backslash\cW_I \mid \Ht(b) < X\} = O_{\eps}\lp\frac{X^{\dim V+\eps}}{\sqrt{M}}\rp + O(X^{\dim V - \delta}).
\]
\end{theorem}

\subsection{A squarefree sieve}
\label{subs:sieve}
Theorem \ref{theo:main} follows from Theorem \ref{theorem:tail} by using
a squarefree sieve. In fact, we will prove a slightly more general result
about families in $\Sigma \subset B(\oO_F)$ defined by infinitely many
congruence conditions.

Let $\kappa$ be a positive integer. We will say that $\cS \subset \Sigma$
is \emph{$\kappa$-acceptable} if $\cS = \cap_{\pp \text{ finite}} \cS_{\pp}$,
where $\cS_{\pp} \subset \Sigma_{\pp} \subset B(\oO_{\pp})$ satisfy the following:
\begin{itemize}
\item $\cS_{\pp}$ is defined by congruence conditions modulo $\pp^{\kappa}$.
\item For all sufficiently large primes $\pp$, the set $\cS_{\pp}$ contains
all elements with $\pp^2$ not dividing $\Delta(b)$.
\end{itemize}
For any subset $A \subset \Sigma$, we denote by $N(A,X)$ the number of
elements of $A$ having height less than $X$. For any prime $\pp$ and any
subset $A_{\pp} \subset \Sigma_{\pp}$, we denote by $\rho(A_{\pp})$ the density
of $A_{\pp}$ inside $\Sigma_{\pp}$.
\begin{theorem}
Let $\cS = \cap_{\pp} \cS_{\pp}$ be a $\kappa$-acceptable subset of $\Sigma$. Then,
there exists a constant $\delta > 0$ such that
\[
N(\cS,X) = \prod_{\pp}\rho(\cS_{\pp}) N(\Sigma,X) + O(X^{\dim V - \delta}).
\]
\end{theorem}
\begin{proof}
Recall that $B = \Spec \oO_F[p_{d_1},\dots,p_{d_k}]$. For an element 
$b \in \Sigma$, we have that $|p_{d_i}(b)| \ll X^{d_i}$, where by
Table \ref{table:explicit} we see that $d_i \geq 2$ for all $i$. For
$I$ a squarefree ideal coprime to $N_{bad}$, we define a family
$\cS^I = \cap_{\pp} \cS_{\pp}^I$ as follows:
\begin{itemize}
\item If $\pp \mid N_{bad}$, then $\cS_{\pp}^I = \cS_{\pp}$.
\item If $\pp \mid I$, then $\cS_{\pp}^I = \Sigma_{\pp} \setminus \cS_{\pp}$.
\item Otherwise, set $\cS_{\pp}^I = \Sigma_{\pp}$.
\end{itemize}
By the inclusion-exclusion principle, we have
\[
N(\cS,X) = \sum_{\substack{I \text{ squarefree}\\(I,N_{bad}) = 1}} \mu(I)N(\cS^I,X),
\]
where $\mu$ is the Möbius function for the ideals of $\oO_F$. By the
Chinese Remainder Theorem, we can estimate
\[
N(\cS^I,X) = \prod_{\pp \mid N_{bad}}\rho(\cS_{\pp}) \prod_{\pp\mid I} (1-\rho(\cS_{\pp})) N(\Sigma,X) + O((NI)^{\kappa} X^{\dim V -2})
\]
From Theorem \ref{theo:fullTail}, we also know that
\[
\sum_{\substack{I \text{ squarefree}\\(I,N_{bad}) = 1\\NI > M}} \mu(I)N(\cS^I,X) = O_{\eps}\lp\frac{X^{\dim V+\eps}}{M}\rp + O(X^{\dim V - \delta}).
\]
Thus, we get that
\begin{align*}
N(\cS^I,X) &= \prod_{\pp \mid N_{bad}}\rho(\cS_{\pp})N(\Sigma,X)\sum_{\substack{I \text{ squarefree}\\(I,N_{bad}) = 1\\NI \leq M}} \mu(I)\prod_{\pp\mid I} (1-\rho(\cS_{\pp})) + O_{\eps}\lp\frac{X^{\dim V+\eps}}{\sqrt{M}} + X^{\dim V - \delta} + M^{\kappa+1}X^{\dim V -2}\rp\\
&= \prod_{\pp}\rho(\cS_{\pp})N(\Sigma,X) + O_{\eps}\lp \frac{X^{\dim V}}{M} + \frac{X^{\dim V+\eps}}{\sqrt{M}} + X^{\dim V - \delta} + M^{\kappa+1}X^{\dim V -2} \rp.
\end{align*}
Optimising, we choose $M = X^{4/(2\kappa + 3)}$ and we conclude the proof.
\end{proof}

\appendix
\section{Counting irreducible orbits}
In the proof of Proposition \ref{prop:Selberg1}, we need a power-saving asymptotic
for the number of orbits in the main body. As in Section \ref{section:Counting},
we let $\Gamma$ be an arithmetic subgroup of $G(F)$, and we let $V_0$
be a lattice of $V(F)$ which is commensurable with $V(\oO_F)$ and
$\Gamma$-stable. We denote by $N^*(\Gamma,V_0,X)$ the number of 
$\Gamma$-orbits in $V_0$ with invariants in $\Lambda\Sigma$ having height at most $X$
(recall that $\Lambda$ is the image of the natural embedding of $\R_{> 0}$
inside $F_{\infty}$). 
Then, we prove the following result:
\begin{theorem}
\label{theo:irred}
There exist positive constants $C,\delta$ such that
\[
N^*(\Gamma,V_0,X) = CX^{\dim V} + O(X^{\dim V - \delta}).
\]
The constant $C$ depends only on $\Gamma$ and $V_0$, while $\delta > 0$
can be chosen independently of $\Gamma$ and $V_0$.
\end{theorem}
\begin{proof}
Following the notations in Section \ref{subs:reductions}, we get that
\[
N^*(\Gamma,V_0,X) = \frac{1}{r_i \vol(G_0)} \sum_{j=1}^m\int_{g \in \alpha_j \cS_j} \#\{v \in V_0 \cap g \cB_X \mid v_{0,j} \neq 0\}dg + O(X^{\dim V - \delta}).
\]
It suffices to estimate the integral for each $j$ separately: for simplicity, set $\alpha_j =1$
and denote the highest weight as $v_0$.
We then have that $|v_0| \geq C_0$ for some constant $C_0$. Let us denote
$\cE = \{v \in g\cB_X \mid |v_0| \geq C_0\}$. There is a constant $J \geq 0$
such that any element in $\cB_1$ satisfies $|v_0| \leq J$, so if $\cE$
is non empty for some choice of $g = ntk$, we must have that $X|a_0(t)| \geq C_0/J$,
where $a_0$ denotes the corresponding weight of $v_0$.

By Davenport's lemma (i.e. Proposition \ref{prop:Davenport}), we can approximate the number of 
lattice points in $\cE$ by some constant times $\vol(\cE)$, with an
error term corresponding to the volume of the lower-dimensional projections. The biggest
volume of a lower-dimensional projection will correspond to setting one
of the real coordinates of $v_0$ to be zero, and this volume can always
be computed to be $X^{\dim V - \delta}$ for a suitable $\delta > 0$.

Then, given that $\cE = g\cB_X \setminus (g\cB_X \setminus \cE)$ it remains to deal with 
\[
\int_{\substack{g = ntk \in \omega_j A_c'K_1\\ X|a_0(t)| \geq C_0/J}} (\vol(g\cB_X) - \vol(g\cB_X \setminus \cE))dg
\]
In the first summand, we note that $\vol(g\cB_X)$ is independent of $g$,
and that $\vol(\cB_X) = X^{\dim V} \vol(\cB_1)$, thus obtaining the main
term. For the second summand, denote $\cE' = g\cB_X \setminus \cE$.
Any element in $v \in g\cB_X \setminus \cE$ 
must have $|v_0| \leq C_0$. Because of how the integration domain is
set up, the values of $v_0$ fall in a compact region $\Omega$ of $F_{\infty}$.
For a given value of $v_0 \in \Omega$, let $\cE'(v_0)$ denote the set
of elements in $\cE'$ with the given value of $v_0$. Then,
\[
\vol(\cE') = \int_{v_0 \in \Omega} \vol(\cE'(v_0)),
\]
and each of the volumes of $\vol(\cE'(v_0))$ can be computed to be
$O(X^{\dim V - \delta})$, for some $\delta > 0$.
\end{proof}
The proof of Theorem \ref{theo:irred} immediately implies the following:
\begin{cor}
\label{cor:irred}
Let $L$ be a translate of $IV_0$ for some ideal $I \subset \oO_F$. Then,
\[
N^*(\Gamma,L,X) = (NI)^{-\dim V} C X^{\dim V} + O((NI)^{-\dim V + \delta} X^{\dim V - \delta}).
\]
Here, $C$ and $\delta$ are as in Theorem \ref{theo:irred}, and the implied
constant is independent of the ideal $I$.
\end{cor}
This result is what we need to apply the Selberg sieve in the proof of
Proposition \ref{prop:Selberg1}.

\printbibliography
\end{document}